\documentclass[reqno,12pt]{amsart}
\usepackage{amsmath,amsthm,amssymb,amsfonts,amscd}
\usepackage{xypic}
\usepackage{appendix}

\usepackage{leftidx}
\def\tr#1{\mathord{\mathopen{{\vphantom{#1}}^t}#1}} 

\theoremstyle{plain} 
	\newtheorem{thm}{Theorem}[section]
	\newtheorem*{thm*}{Theorem}
	
	\newtheorem{lem}[thm]{Lemma}
	
	\newtheorem{prop}[thm]{Proposition}
	
	\newtheorem*{conj*}{Conjecture}
	
\theoremstyle{definition}
	\newtheorem{defn}[thm]{\rm{Definition}}

\theoremstyle{remark}
	\newtheorem{rmk}[thm]{\rm{Remark}}
	
	\newtheorem*{pf}{\rm{Proof}}
\numberwithin{equation}{section}
\def\CC{{\mathbb C}}

\def\HH{{\mathbb H}}

\def\PP{{\mathbb P}}
\def\QQ{{\mathbb Q}}
\def\RR{{\mathbb R}}
\def\ZZ{{\mathbb Z}}

\def\F{{\mathcal F}}

\def\I{{\mathcal I}}

\def\O{{\mathcal O}}

\def\p{\partial }
\def\Oo{{\O}}

\newcommand{\bp}{\begin{pmatrix}}
\newcommand{\ep}{\end{pmatrix}}

\newcommand{\q}{(0)}

\numberwithin{equation}{section}

\newcounter{CounterEQUlabel}
\newcommand{\EQUlabel}[1]{\label{#1}
	\ifcase \theCounterEQUlabel
		\relax
	\or
		\hspace{1em}\mbox{\tiny$\langle$\rmfamily#1$\rangle$}
		\index{zzz#1@#1}
	\fi }	
	\newcounter{CounterEQUref}
	\newcounter{CounterEQUpageref}
	\newcommand{\EQUref}[1]{
		\ifcase \theCounterEQUref     \relax   \or {\tiny[#1]}\,\fi
		\ifcase \theCounterEQUpageref (\ref{#1}) \or (\ref{#1}\,(p.\pageref{#1})) \fi}

\begin{document}
\title{Frobenius structures and characters of affine Lie algebras}
\date{\today}
\author{Ikuo Satake}
\address{Faculty of Education, Kagawa University, 
1-1 Saiwai-cho Takamatsu Kagawa, 760-8522, Japan}
\email{satakeikuo@gmail.com}
\begin{abstract}
The explicit description of the Frobenius structure for 
the elliptic root system of type $D_4^{(1,1)}$ 
in terms of the characters of an affine Lie algebra 
of type $D_4^{(1)}$ is given. 
\end{abstract}
\maketitle
\section{Introduction}

As a generalization of root systems, 
the elliptic root systems are defined in 
Saito \cite{extendedI}. 
The Frobenius structure (also called ``the flat structure") 
for the elliptic root system 
is introduced \cite{extendedII, handai}. 
One of the purposes of introducing this structure is 
to construct the automorphic functions (cf. Saito \cite{period}), 
where the explicit description of the Frobenius structure 
is important. 
For the elliptic root systems of types 
$G_2^{(1,1)}$ \cite{Kokyuroku, Bertola}, 
$D_4^{(1,1)}$ \cite{D4} and 
$E_6^{(1,1)}$ \cite{E6}, 
the results of the explicit descriptions of the Frobenius structures are given. 

In Satake-Takahashi \cite{ST}, the description of 
the Frobenius structure for the elliptic root system of type $D_4^{(1,1)}$ in \cite{D4} 
is compared with the orbifold Gromov-Witten invariants 
for the elliptic orbifold of type $(2,2,2,2)$. 
However, the proofs of the results of \cite{D4} are omitted. 
In this paper, we give the proofs of the results of \cite{D4}. 

The contents of this paper are as follows. 
In Section 2, we recall the elliptic root system \cite{extendedII}. 
Also, the normalized characters of affine Lie algebras are introduced 
in this setting. 
In Section 3, we review the construction of the Frobenius structure 
for the elliptic root system \cite{extendedII, handai}. 
In Section 4, the main theorem (Theorem \ref{400.001}) 
which gives the explicit description 
of a flat generator system and the potential of the Frobenius structure 
for the elliptic root system of type $D_4^{(1,1)}$ is given. 

In Section 5, the proof of the main theorem is given. 
The relation of the Frobenius structure and 
the normalized characters of affine Lie algebras 
are, as formula, very simple. 
However, the conceptual relation among these structures 
is not yet known. 
Thereby, we first construct ``Jacobi forms" (\ref{520.001}) 
from the normalized characters. 
Then the relation between these Jacobi forms and the Frobenius structure 
is studied 
in Proposition \ref{170831.31} 
using differential relations satisfied by Jacobi forms 
which are given in Section 6. 
In Proposition \ref{170831.33}, the relations 
of the normalized characters of affine Lie algebras 
and the Frobenius structure are studied, where we need 
the analyses of the respective solutions of both the 
Kaneko-Zagier equations and the Halphen's equations which are given in Section 7. 

For the Jacobi forms, the integers called ``weight" and ``index" are 
defined and their ``initial terms" are also defined. 
However, these notions do not fit the differential operators 
which appear in the construction of the Frobenius structure. 
Thus, in Section 6 we formulate the (lifted) Jacobi forms and the differential operators 
on a $1$-dimensional higher space where they fit together well. 
Then we obtain results such as ``the determinant construction" of the Jacobi forms 
from normalized characters in Section \ref{622.002} and 
a description of the intersection form of the Frobenius structure. 
Then we give a reduction to the $1$-dimensional lower space discussed in Section 5. 
This is the reason for separating Section 5 and Section 6. 


In Section 7, the Kaneko-Zagier equations and 
the Halphen's equations satisfied by modular forms are studied. 
In Section 8, the duality 
of the solutions of 
the Kaneko-Zagier equations is formulated. 
These results allow for the description of the Frobenius structure. 
In Appendix A, we give a lemma on the connection 
used to construct flat sections in Proposition \ref{170831.32}. 

\textrm{Acknowledgements:}
The author would like to thank the referee 
for careful reading of the manuscript. 
The author is supported by JSPS KAKENHI Grant Number 
JP17K18781, JP15K13234.

%

\section{$W$-invariants for elliptic root systems}
We recall the elliptic root systems (cf. Saito \cite{extendedII}). 
Some parts are changed in order to fit the notations of Kac \cite{Kac}. 

\subsection{Elliptic root systems}

In this subsection, we define the elliptic root systems (cf. \cite{extendedII}). 

Let $l$ be a positive integer. 
Let $F$ be a real vector space of rank $l+2$ with a positive semi-definite 
symmetric bilinear form 
$I:F \times F \to \RR$, whose radical 
$\mathrm{rad}I:=\{x \in F\,|\,I(x,y)=0,\forall y \in F\}$ 
is a vector space of rank $2$. 
For a non-isotropic vector $\alpha \in F$ 
(i.e. $I(\alpha,\alpha)\neq 0$), 
we put $\alpha^{\vee}:=2\alpha/I(\alpha,\alpha)\in F$. 
The reflection $w_{\alpha}$ with respect to $\alpha$ is 
defined by 
$
w_{\alpha}(u):=u-I(u,\alpha^{\vee})\alpha
\quad
(\forall u \in F)
$. 
\begin{defn}
A set $R$ of non-isotropic elements 
of $F$ is an elliptic root system belonging to $(F,I)$ if it satisfies 
the axioms (i)-(iv):
\begin{enumerate}
\item The additive group generated by $R$ in $F$, 
denoted by $Q(R)$, 
is a full sub-lattice of $F$.
\item $I(\alpha,\beta^{\vee}) \in \ZZ$ for $\alpha,\beta \in R$. 
\item $w_{\alpha}(R)=R$ for $\forall \alpha \in R$. 
\item If $R=R_1 \cup R_2$, with $R_1 \perp R_2$, then either 
$R_1$ or $R_2$ is void. 
\end{enumerate}
\end{defn}
We have $Q(R) \cap \mathrm{rad}I \simeq\ZZ^2$. 
We call a 1-dimensional vector space $G \subset \mathrm{rad}I$ 
satisfying $G \cap Q(R) \simeq \ZZ$, a marking. 
We fix $a,\delta \in F$ s.t. 
$G \cap Q(R)=\ZZ a$ and $Q(R) \cap \mathrm{rad}I=\ZZ a \oplus \ZZ \delta$. 

The isomorphism classes of the elliptic root systems 
with markings are classified in \cite{extendedI}. 
Hereafter we restrict ourselves to the elliptic root systems 
of type $X_l^{(1,1)}$ $(X_l=A_l,B_l,\cdots,G_2)$ in the notation of \cite{extendedI}. 

For these elliptic root systems, an $\mathbb{R}$-basis of $F$ is given as follows. 
We have a natural projection $p:F \to F/\mathrm{rad}I$. 
Put $R_{\mathrm{cl}}:=p(R)$. 
Then, $R_{\mathrm{cl}} \subset F/\mathrm{rad}I$ with an $\RR$-bilinear form 
on $F/\mathrm{rad}I$ induced by $I$ gives a finite root system. 
We take $\alpha_1,\cdots,\alpha_l \in R$ s.t. 
$p(\alpha_1),\cdots,p(\alpha_l)$ give simple roots, where 
we follow the enumeration of the vertices of the Dynkin diagram of Kac \cite[p.53]{Kac}. 
We put $F_{\mathrm{cl}}:=\bigoplus_{i=1}^l \RR\alpha_i$. 
Then we have $F=F_{\mathrm{cl}}\oplus \RR \delta \oplus \RR a$. 

We fix some notations. 
There exists a unique $\theta \in F_{\mathrm{cl}} \cap R$ which satisfies 
the condition that $p(\theta)$ is the highest root of the finite root 
system $p(R)\subset F/\mathrm{rad}I$. 
Hereafter we assume that $I(\theta,\theta)=2$ which is satisfied by a constant 
multiplication of $I$ if necessary. 
For $\alpha \in R$, we put $\alpha^{\vee}=2\alpha/I(\alpha,\alpha)$. 
For $i=1,\cdots,l$, we define $a_i^{\vee}\in \RR$ by the coefficients of 
$\theta=\sum_{i=1}^l a_i^{\vee}\alpha^{\vee}_i$. Also we put $a_0^{\vee}=1$. 
These $a_i^{\vee}$ $(i=0,\cdots,l)$ are known to be positive integers and called 
colabels. 
We put $\alpha_0=\delta-\theta$. 
We have $\alpha_0 \in R$. 
Then we have $\delta=\sum_{i=0}^l a_i^{\vee}\alpha^{\vee}_i$. 

\subsection{Hyperbolic extension and fundamental weights}

We introduce a hyperbolic extension $(\tilde{F},\tilde{I})$ of $(F,I)$, 
i.e. $\tilde{F}$ is a $(l+3)$-dimensional $\RR$-vector space of 
which contain $F$ as a subspace 
and $\tilde{I}$ is a symmetric $\RR$-bilinear form 
on $\tilde{F}$ which satisfies $\tilde{I}|_F=I$ and $\mathrm{rad}\tilde{I}=\RR a$. 
It is unique up to isomorphism. 

We fix some notations. 
We define $\Lambda_0 \in \tilde{F}$ which satisfies 
$\tilde{I}(\Lambda_0,\alpha_i^{\vee})=\delta_{0,i}$ for $i=0,\cdots,l$ 
and $\tilde{I}(\Lambda_0,\Lambda_0)=0$. 
Then we have a decomposition $\tilde{F}=F\oplus \RR\Lambda_0$. 
For $i=1,\cdots,l$, we define $\omega_i \in F_{\mathrm{cl}}$ 
by $I(\omega_i,\alpha_j^{\vee})=\delta_{ij}$ for $j=1,\cdots,l$. 
We put $\Lambda_i:=a_i^{\vee}\Lambda_0+\omega_i \in \tilde{F}$ 
for $i=1,\cdots,l$. 
We call $\Lambda_0,\cdots,\Lambda_l$, the fundamental weights. 
We remark that $I(\Lambda_i,\alpha_j^{\vee})=\delta_{ij}$ for 
$i,j=0,\cdots,l$.  
We denote the sum of fundamental weights $\sum_{i=0}^l \Lambda_i$ by $\rho$. 

\subsection{The elliptic Weyl group and the affine Weyl group}
We define the elliptic Weyl group. 
We first define 
$
O(\tilde{F},\mathrm{rad}I)
:=
\{g \in GL(\tilde{F})\,|\,
\tilde{I}(gx,gy)=\tilde{I}(x,y)\ 
\forall x,y \in \tilde{F},\ 
g|_{\mathrm{rad}I}=id. \}$. 
For $\alpha \in R$, we define a reflection $\tilde{w}_{\alpha} \in 
O(\tilde{F},\mathrm{rad}I)$ by 
$\tilde{w}_{\alpha}(u)=u-\tilde{I}(u,\alpha^{\vee})\alpha$ for $u \in \tilde{F}$.  
We define the elliptic Weyl group $W$ by the 
group generated by $\tilde{w}_{\alpha}\ (\alpha \in R)$. 
We have a homomorphism:
$pr_*:O(\tilde{F},\mathrm{rad}I) \to GL(\tilde{F}/\RR a)$
induced by $pr:\tilde{F} \to \tilde{F}/\RR a$. 
Also we have a homomorphism $\oplus_{i=0}^l\ZZ\alpha_i 
\to W $ which is defined by $\sum k_i \alpha_i^{\vee} \mapsto 
(u \mapsto u-a\sum \tilde{I}(\sum k_i \alpha_i^{\vee},u))$. 
Then we have the following exact sequence (cf. Saito-Takebayashi \cite{Saito-Takebayashi}):
\begin{equation}\label{230.001}
0 \to \bigoplus_{i=0}^l\ZZ\alpha_i^{\vee} 
\to W \to pr_*(W)\to 1.
\end{equation}

We also define the affine Weyl group $W_{\mathrm{af}}$ by 
the group generated by $\tilde{w}_{\alpha_i}\ (i=0,\cdots,l)$. 
Then $W_{\mathrm{af}}$ gives a splitting of the exact sequence (\ref{230.001}) 
and we have a semi-direct product decomposition:
\begin{equation}
W=W_{\mathrm{af}} \ltimes (\bigoplus_{i=0}^l\ZZ\alpha_i^{\vee} ).\label{230.002}
\end{equation}
We have a group homomorphism $\varepsilon:W_{\mathrm{af}} \to \{\pm 1\}$ 
defined by the parity of number of reflections 
which define $w \in W_{\mathrm{af}}$.

\subsection{Domains and a bilinear form}\label{240.001}

We define two domains:
\begin{eqnarray*}
Y&&:=\{x \in \mathrm{Hom}_{\RR}(\tilde{F},\CC)\,|\,\langle a,x\rangle=-2\pi\sqrt{-1},\ 
		\mathrm{Re}\langle \delta,x\rangle>0\},\\
H&&:=\{x \in \mathrm{Hom}_{\RR}(\mathrm{rad}I,\CC)\,|\,\langle a,x\rangle=-2\pi\sqrt{-1},\ 
		\mathrm{Re}\langle \delta,x\rangle>0\}.
\end{eqnarray*}
We have a natural morphism $\pi:Y\to H$. 
Hereafter we shall identify 
$H$ and $\mathbb{H}:=\{z\in \CC\,|\,\mathrm{Im}z>0\}$ 
by the function 
$\delta/a:H \to \mathbb{H}$. 
We define the left action of $g \in W$ on $Y$ by 
$\langle g\cdot x,\gamma \rangle=\langle x,g^{-1} \cdot \gamma \rangle$. 
for $x \in Y$ and $\gamma \in \tilde{F}$. 
We remark that the domain $Y$ and the action of $W$ and $W_{\mathrm{af}}$ on $Y$ 
is naturally identified with Kac \cite[p.225]{Kac}. 

We shall regard 
$\delta,\alpha_1,\cdots,\alpha_l,\Lambda_0$ as a coordinate system 
of $Y$. 

We define a vector field $E$ by 
$$
E:=\frac{\p}{\p \Lambda_0}:\Omega(Y) \to F(Y), 
$$
where $F(Y)$ is a space of holomorphic functions on $Y$ 
and $\Omega(Y)$ is a space of holomorphic 1-forms on $Y$. 
Put $m_{\mathrm{max}}:=\mathrm{max}\{a_0^{\vee},\cdots,a_l^{\vee}\}$. Then 
we define $E_{\mathrm{norm}}:=\frac{1}{m_{\mathrm{max}}}E$. 
\begin{rmk}
$E_{\mathrm{norm}}$ is used to define Euler field 
for the Frobenius structure. 
$E$ is used to count the degree of the $W$-invariants. 
\end{rmk}

We define a symmetric bilinear form 
$I^*:\Omega(Y) \times \Omega(Y) \to F(Y)$ 
and its Laplacian $D$ and $\overline{D}$ by 
$$
I^*:=-
\left[
\frac{\p}{\p \delta}\otimes 
\frac{\p}{\p \Lambda_0}
+ 
\frac{\p}{\p \Lambda_0}
\otimes 
\frac{\p}{\p \delta}
+\sum_{i,j=1}^l
\tilde{I}(\alpha_i,\alpha_j)
\frac{\p}{\p \alpha_i}
\otimes 
\frac{\p}{\p \alpha_j}
\right],
$$
$$
\overline{D}:=
\sum_{i,j=1}^l
\tilde{I}(\alpha_i,\alpha_j)
\frac{\p^2}{\p \alpha_i\p \alpha_j}
,
$$
$$
D:=
-
\left[
2
\frac{\p^2}{\p \delta \p \Lambda_0}
+\overline{D}
\right].
$$

We define a map 
$$
Y \to \mathbb{H}, \quad 
h \mapsto \frac{\delta(h)}{a(h)}=\frac{\delta(h)}{-2\pi\sqrt{-1}}, 
$$
and denote it by $\tau$. 
We also define maps 
\begin{align*}
&Y \to \mathrm{Hom}_{\RR}(F_{\mathrm{cl}},\CC), \quad
h \mapsto \left(\gamma \mapsto -\frac{\gamma(h)}{a(h)}
=\frac{\gamma(h)}{2\pi\sqrt{-1}}\right),\\
&Y \to \CC, \quad
h \mapsto -\frac{\Lambda_0(h)}{a(h)}=\frac{\Lambda_0(h)}{2\pi\sqrt{-1}}. 
\end{align*}
Thus we have 
\begin{equation}
Y \simeq \mathbb{H}\times \mathfrak{h}_{\CC} \times \CC,\label{240.003}
\end{equation}
where we denote $\mathrm{Hom}_{\RR}(F_{\mathrm{cl}},\CC)$ by 
$\mathfrak{h}_{\CC}$ which is isomorphic to the 
complexified Cartan subalgebra of the finite dimensional Lie algebra 
corresponding to the root system $R_{\mathrm{cl}}$. 

We remark that the symmetric bilinear form $I^*$ induces 
a morphism $\Omega(Y) \to \Theta(Y)$ 
where $\Theta(Y)$ is a space of $\CC$-derivations of $S^W$ 
and we also denote it by $I^*$. Then we have 
\begin{equation}\label{240.002}
I^*(d(2\pi\sqrt{-1}\tau))=E. 
\end{equation}
\subsection{The algebra of the invariants for the elliptic Weyl group}
In this subsection, we introduce the algebra of the invariants 
for the elliptic Weyl group. 

We put $F(\mathbb{H}):=\{f:\mathbb{H} \to \CC:holomorphic\}$. 
For $m \in\ZZ$, we put 
$F_m(Y):=\{f\in F(Y)\,|\,Ef=mf\}$. 
By $\pi^*:F(\mathbb{H}) \to F_0(Y)$, induced 
by $\pi:Y \to H$, $F_m(Y)$ is a $F(\mathbb{H})$-module. 

For $m \in \ZZ$ , we put 
$S_m^W:=\{f\in F_m(Y)\,|\,f(g\cdot z)=f(z)\ \forall g \in W\}$, 
$S^W:=\bigoplus_{m \in \ZZ_{\geq 0}}S^W_{m}$. 
$S^W$ is a graded algebra and we see that $S^W_0$ is just 
$\pi^*(F(\mathbb{H}))$. 

\begin{thm}[\cite{Chevalley3, Chevalley4, Chevalley2, Chevalley1, Chevalley5}]
There exist $l+1$ homogeneous elements $P_i \in S^W_{a_i^{\vee}}$ 
$(0 \leq i \leq l)$ 
such that $S^W$ is a polynomial algebra 
of $P_0,\cdots,P_l$ over $F(\mathbb{H})$. 
\end{thm}

\subsection{$S^W$-modules and $S^W$-bilinear forms}\label{170930}
We introduce $S^W$ modules and tensors on them.  

We define $\Omega_{S^W}$ as a set of 1-forms for the algebra $S^W$. 
We define $\Theta_{S^W}$ as $\CC$-derivations of the algebra $S^W$. 
We have 
$$
\Theta_{S^W}=\bigoplus_{i=-1}^l S^W \frac{\p}{\p P_i}, 
\ 
\Omega_{S^W}=\bigoplus_{i=-1}^l S^W dP_i, \ 
P_{-1}:=\tau.
$$
Thus $\Omega_{S^W},\ \Theta_{S^W}$ 
are free graded $S^W$-modules of rank $l+2$. 

We introduce a filtration on $\Omega_{S^W}$ by the degree of the 
$S^W$-free generators:
\begin{equation}\label{250.001}
\{0\}=F_{-1}\subset S^W\frac{dq}{q}=F_0 \subset F_1 
\subset \cdots \subset F_{m_{\mathrm{max}}}
=\Omega_{S^W},
\end{equation}
where 
$$
F_k:=\bigoplus_{\deg P_i\leq k}S^W dP_i \quad(-1 \leq k \leq m_{\mathrm{max}}), 
$$
which do not depend on the choice of $P_i$. 

We denote by $q$ the holomorphic function $e^{-\delta}=e^{2\pi\sqrt{-1}\tau}$ 
on $\HH$. 
Then we have $\frac{dq}{q}=2\pi\sqrt{-1}dP_{-1}=2\pi\sqrt{-1}d\tau$, 
$q\frac{d}{dq}=\frac{1}{2\pi\sqrt{-1}}\frac{\p}{\p P_{-1}}
=\frac{1}{2\pi\sqrt{-1}}\frac{\p}{\p \tau}$. 
For $f \in F(\HH)$, we denote $q\frac{d}{dq}f$ by $f'$. 

By the natural morphism $\Omega_{S^W} \to \Omega(Y)$, 
$E,I^*,D$ defined in Section \ref{240.001} restricts: 
\begin{eqnarray}
E:&&\Omega_{S^W}\to S^W,\\
E_{\mathrm{norm}}:&&\Omega_{S^W}\to S^W,\\
I^*:&&\Omega_{S^W}\times \Omega_{S^W}\to S^W,\\
D:&&\Omega_{S^W}\to S^W,
\end{eqnarray}
since $Lie_{E}E=0$, $Lie_{E}I^*=0$. 
We call $E_{\mathrm{norm}}$ the Euler field. 

For $f \in F(\mathbb{H})$, $F \in F_m(Y)$ $(m \in \ZZ)$, 
we have 
\begin{eqnarray}
D(fF)&=&2(q\frac{d}{dq}f)(EF)+fD(F),\label{250.002}\\
I^*(df,dg)&=&\frac{1}{2}[D(fg)-fD(g)-gD(f)].\label{250.003}
\end{eqnarray}

\subsection{The normalized characters}\label{260.002}
We introduce the normalized characters $\chi_{\Lambda_i}$ $(i=1,\cdots,l)$ 
motivated by Kac \cite[p.226]{Kac}. 
\begin{defn} 
For $i=0,1,\cdots,l$, 
we define holomorphic functions on $Y$: 
\begin{eqnarray}
ch_{\Lambda_i}&:=&(\sum_{w \in W_{\mathrm{af}}} \varepsilon(w)e^{w(\rho+\Lambda_i)})
/(\sum_{w \in W_{\mathrm{af}}} \varepsilon(w)e^{w(\rho)}), \\
\chi_{\Lambda_i}&:=&e^{-m_{\Lambda_i}\delta}ch_{\Lambda_i}
=q^{m_{\lambda_i}}ch_{\Lambda_i}, 
\end{eqnarray}
where $m_{\Lambda_i}:=\frac{|\Lambda_i+\rho|^2}
{2\tilde{I}(\Lambda_i+\rho,\delta)}
-\frac{|\rho|^2}{2\tilde{I}(\rho,\delta)}$. 
\end{defn}
The function $\chi_{\Lambda_i}$ is 
invariant by the action of $\bigoplus_{i=0}^l{\mathbb Z}\alpha_i^{\vee}$ 
and $W_{\mathrm{af}}$ of the subgroups $W$. 
Then by (\ref{230.002}), we have 
\begin{prop}\label{260.001}
For $i=0,1,\cdots,l$, $\chi_{\Lambda_i}\in S^W_{a_i^{\vee}}$.
\end{prop}

\subsection{Some modular forms}
We put 
\begin{eqnarray}
E_2&:=&1-24\sum_{n=1}^{\infty}\sigma_1(n)q^n , \label{270.001}\\
\eta&:=&q^{1/24}\prod_{n=1}^{\infty}(1-q^n) , \label{270.002}
\end{eqnarray}
where $\sigma_k(n):=\sum_{d|n}d^k$. 
$E_2,\,\eta$ are holomorphic functions on $\HH$ 
and we have $E_2=24(\eta'/\eta)$. 

We define $ch_{\Lambda_i}^{\q}$ and $\chi_{\Lambda_i}^{\q} \in F(\HH)$ 
by $ch_{\Lambda_i}(\tau,0,0)$ and $\chi_{\Lambda_i}(\tau,0,0)$ 
respectively, where $ch_{\Lambda_i}(\tau,z,t),\,\chi_{\Lambda_i}(\tau,z,t)$ 
are defined through the isomorphism 
$(\tau,z,t) \in \HH \times \mathfrak{h}_{\CC} \times \CC\simeq Y$. 
%

\section{Frobenius structure }

The following theorems are proved in \cite{handai} 
in the context of Frobenius manifolds. 
\begin{thm}\label{300.001}
Assume that $X_l^{(1,1)}=D_4^{(1,1)},E_6^{(1,1)},E_7^{(1,1)},
E_8^{(1,1)},F_4^{(1,1)},G_2^{(1,1)}$. 
Then there exist the holomorphic metric 
$J:\Theta_{S^W} \otimes_{S^W} \Theta_{S^W} \to {S^W}$ and 
multiplication 
$\circ:\Theta_{S^W} \otimes_{{S^W}}\Theta_{S^W} \to \Theta_{S^W}$ 
on $\Theta_{S^W}$, global unit field $e:\Omega_{S^W} \to {S^W}$, 
satisfying the following conditions:
\begin{enumerate}
\item the metric is invariant under the multiplication, 
i.e. $J(X\circ Y,Z)=J(X,Y \circ Z)$ 
for the vector field $X,Y,Z:\Omega_{S^W} \to {S^W}$, 
\item $($potentiality$)$ the 
$(3,1)$-tensor $\nabla \circ$ is symmetric
(here $\nabla$ is the Levi-Civita connection of the metric), 
i.e. 
$\nabla_X(Y \circ Z)
-Y\circ \nabla_X (Z)
-\nabla_Y(X \circ Z)
+X\circ \nabla_Y (Z)
-[X,Y]\circ Z=0$, 
for the vector field $X,Y,Z:\Omega_{S^W} \to {S^W}$, 
\item the metric $J$ is flat, 
\item $e$ is a unit field and it is flat, i.e. $\nabla e=0$,
\item the Euler field $E_{\mathrm{norm}}$ satisfies 
$Lie_{E_{\mathrm{norm}}}(\circ)=1 \cdot \circ$, and 
$Lie_{E_{\mathrm{norm}}}(J)=1 \cdot J$, 
\item the intersection form coincides with $I^*$:
$J(E_{\mathrm{norm}},J^*(\omega)\circ J^*(\omega'))=I^*(\omega,\omega')$ 
for 1-forms $\omega,\omega' \in \Omega_{S^W}$ and 
$J^*:\Omega_{S^W} \to \Theta_{S^W}$ is the isomorphism induced by $J$. 
\end{enumerate}
\end{thm}

\begin{thm}\label{300.002}
For $c \in \CC^*$, $(c^{-1}\circ,ce,c^{-1}J)$ 
also satisfy the conditions Theorem \ref{300.001}. 
If $(\circ',e',J')$ satisfy conditions Theorem \ref{300.001}, 
then there exists $c \in \CC^*$ such that 
$(\circ',e',J')=(c^{-1}\circ,ce,c^{-1}J)$. 
In particular the Levi-Civita connection $\nabla$ is unique 
and the vector space $V:=\CC e$ does not depend on the choice of $e$. 
\end{thm}

We introduce the potential $\F \in S^W$ of the Frobenius structure 
by the condition that 
$E_{\mathrm{norm}}\F=2\F$, $J(X \circ Y,Z)=XYZ \F$ 
for $X,Y,Z \in \Theta_{S^W}$ with $\nabla X=\nabla Y=\nabla Z=0$. 

We introduce the flat generator system as follows. 
Put $t_{-1}:=2\pi\sqrt{-1}\tau/m_{\mathrm{max}}$ and $d_{-1}=0$. 
Then we have $\nabla dt_{-1}=0$ and 
$E_{\mathrm{norm}}t_{-1}=d_{-1} t_{-1}$. 

Let $t_i \in S^W$ $(i=0,\cdots,l)$ be such that $\nabla dt_i=0$, 
$E_{\mathrm{norm}}t_i=d_it_i$  
and $dt_{-1},dt_0,\cdots,dt_l$ form an $S^W$-free basis of $\Omega_{S^W}$. 
Then we call $t_{-1},t_0,\cdots,t_l$ the flat generator system. 
By Theorem \ref{300.001}, 
we assume that $0<d_0\leq \cdots \leq d_{l-1}<d_l=1$. 

By the flat generator system, 
the potential $\F$ is uniquely determined up to adding a constant multiple of $t_l^2$ 
by the inequality $d_{l-1}<d_l=1$. 

The following Proposition will be used to construct 
the potential $\F$. 
\begin{prop}\label{510.004}
For $i \neq -1$ or $j\neq -1$, we have 
$\p^k\p^i\p^j\F=\p^k[\frac{1}{d_i+d_j}I^*(dt_i,dt_j)]$ for 
$\p^k=J^*(dt^k)$. 
For $k=i=j=-1$, $\p^{k}\p^{i}\p^{j}\F=0$. 
\end{prop}

\begin{pf}
Since $\p^{i}$ is of degree $d_i-1$, the degree 
of $\p^{k}\p^{i}\p^{j}\F$ is $d_k+d_i+d_j-3+2$, we have 
$\p^{k}\p^{i}\p^{j}\F=0$ for $k=i=j=-1$. 
By Theorem \ref{300.001} (vi) and definition of $\F$, 
we have $I^*(dt_i,dt_j)=E_{\mathrm{norm}}\p^i\p^j\F=(d_i+d_j)\p^i\p^j\F$. 
If $i\neq -1$ or $j\neq -1$, then $d_i+d_j \neq 0$. 
Then we have the result. 
\qed
\end{pf}

\section{Main theorem}

\begin{thm}\label{400.001}
For the elliptic root system of type $D_4^{(1,1)}$, the following holds. 
\begin{enumerate}
\item Put 
$
b_{-1}:=\pi\sqrt{-1}\tau,
b_0:=\eta^{-2}\chi_{\Lambda_0},\ 
b_1:=\eta^{-2}\chi_{\Lambda_1},\ 
b_2:=\eta^{-2}\chi_{\Lambda_3},
b_3:=\eta^{-2}\chi_{\Lambda_4},
$
$$
b_4:=\frac{1}{8}\sum_{i=0,1,2,3}
I^*(db_i,db_i)
+\frac{3}{2}(b_0^2+b_1^2+b_2^2+b_3^2)\frac{\eta'}{\eta}.
$$
Then $b_0,\cdots,b_4$ are polynomial generators 
of $S^W$ over $F(\HH) $ and 
$\nabla db_i=0$ $(i=-1,\cdots,4)$ 
with respect to $\nabla$ defined in Theorem $\ref{300.002}$. 
Namely, $b_i$ $(i=-1,\cdots,4)$ are homogeneous flat coordinates. 
\item We take $\frac{\p}{\p b_i}$ $(i=-1,\cdots,4)$ 
$S^W$-free basis of $\Theta_{S^W}$. 
Then $\frac{\p}{\p b_4}\in V$. 
Put $e_0:=\frac{\p}{\p b_4}$. By Theorem $\ref{300.002}$, 
we have the unique Frobenius structure whose unit of the product is $e_0$. 
We denote its holomorphic metric by $J_0$. 
Then we have 
$$
(J_0^*(db_i,db_j))=
\bp 
0&0&0&0&0&1\\
0&2&0&0&0&0\\
0&0&2&0&0&0\\
0&0&0&2&0&0\\
0&0&0&0&2&0\\
1&0&0&0&0&0
\ep. 
$$
The following function 
\begin{eqnarray*}
\F_0&:=&
\frac{1}{2}(\pi \sqrt{-1}\tau)(b_4)^2+\frac{1}{4}b_4(b_0^2+b_1^2+b_2^2+b_3^2)\\
&&+f_0(b_0b_1b_2b_3)\\
&&+\frac{1}{4}f_1(b_0^4+b_1^4+b_2^4+b_3^4)\\
&&+\frac{1}{6}
f_2(b_0^2b_1^2+b_0^2b_2^2+b_0^2b_3^2+b_1^2b_2^2+b_1^2b_3^2+b_2^2b_3^2)
\end{eqnarray*}
gives a potential for this Frobenius structure where 
\begin{eqnarray*}
f_0&:=&\frac{1}{8}\eta^4\chi_{\Lambda_1}^{\q},\\
f_1&:=&-\frac{1}{2}\left(\frac{1}{24}E_2+\frac{1}{24}\eta^4\chi_{\Lambda_0}^{\q}\right),\\
f_2&:=&-\frac{3}{2}\left(\frac{1}{24}E_2-\frac{1}{24}\eta^4\chi_{\Lambda_0}^{\q}\right).
\end{eqnarray*}
\end{enumerate}
\end{thm}

\begin{rmk}
In Satake \cite{D4} and Satake-Takahashi \cite{ST}, the corresponding results are described. 
The relation of 
$\eta^4\chi_{\Lambda_0}^{\q}, \eta^4\chi_{\Lambda_1}^{\q}$ 
and 
$\Theta_{0,1}, \Theta_{\omega_1,1}$ in \cite{ST} 
is 
$\Theta_{0,1}(e^{\pi\sqrt{-1}\tau})
=\eta^4(e^{2\pi\sqrt{-1}\tau})
\chi_{\Lambda_0}^{\q}(e^{2\pi\sqrt{-1}\tau})$, 
$\Theta_{\omega_1,1}(e^{\pi\sqrt{-1}\tau})
=\eta^4(e^{2\pi\sqrt{-1}\tau})
\chi_{\Lambda_1}^{\q}(e^{2\pi\sqrt{-1}\tau})$ 
because of the result of string function in \cite{Chevalley2}. 
We remark that $q$ in \cite{ST} is $e^{\pi \sqrt{-1}\tau}$ 
whereas $q$ is $e^{2\pi\sqrt{-1}\tau}$ in this paper. 
Also the functions $X_i$ in Lemma 2.3 of \cite{ST} 
is the same of the functions $\xi_i$ 
in Section \ref{170831.3} (\ref{641.004}). Thus we may identify 
our results with those in \cite{ST}. 
\end{rmk}

\section{Proof of main theorem}

\subsection{Characterization of the unit field}
In this subsection, we discuss under the setting of Section 3. 
We give a characterization of 
the unit field $e$ up to constant multiple which 
will be used in Proposition \ref{170831.31}. 

We first remind following general results for the Frobenius structure. 
\begin{prop}\label{510.001}
Let $J^*:\Omega_{S^W}\times \Omega_{S^W} \to S^W$ be a dual of $J$. 
Then we have $J^*=Lie_eI^*$. 
Also we have $Lie_eLie_eI^*=0$.  
\end{prop}
\begin{pf}
By Hertling \cite[p.146]{Hertling}, $Lie_eJ=0$, $Lie_e \circ=\circ$, $Lie_e E_{\mathrm{norm}}=e$. 
Then we have the results. 
\qed\end{pf}

We return to the situation of 
the Frobenius structures for elliptic root systems. 
We put 
$V_1:=\{v \in \Theta_{S^W}\,|\,Lie_{E_{\mathrm{norm}}}v=-v,\ \nabla v=0\}$, 
$V_2:=\{v \in \Theta_{S^W}\,|\,Lie_{E_{\mathrm{norm}}}v=-v,\ Lie_vLie_vI^*=0\}$. 
\begin{prop}\label{510.003}
For $V$ in Theorem $\ref{300.002}$, we have $V=V_1=V_2$. 
\end{prop}
\begin{pf}
We see $V \subset V_1$, $V \subset V_2$. 
For the case the Froenius structure for the elliptic root system, 
the degrees of flat generator system are described by $0,\,a_i^{\vee}/m_{\mathrm{max}}$ 
$(0 \leq i \leq l)$ 
and we see that $\mathrm{dim}_{\CC}V_1=1$. 
The fact that $\mathrm{dim}_{\CC}V_2=1$ is shown in Saito \cite[\S 9.3]{extendedII}. 
Thus we have the result. 
\qed\end{pf}

\subsection{Proof of main theorem}
We give a proof of Theorem \ref{400.001} by using 
Proposition \ref{170831.31}, \ref{170831.34}, \ref{170831.32}, \ref{170831.33}. 

We put 
\begin{subequations}\label{520.001}
\begin{align}
s_{-1}&:=\pi\sqrt{-1}\tau,\\
s_0&:=
(-6)\eta^{-4}
\left[(\chi_{\Lambda_0}^{\q})'(\chi_{\Lambda_1}+\chi_{\Lambda_3}+\chi_{\Lambda_4})
-(\chi_{\Lambda_1}^{\q}+\chi_{\Lambda_3}^{\q}+\chi_{\Lambda_4}^{\q})'(\chi_{\Lambda_0})
\right],\\
s_1&:=
\eta^{-4}
[(\chi_{\Lambda_1}+\chi_{\Lambda_3}+\chi_{\Lambda_4})\chi_{\Lambda_0}^{\q}
-\chi_{\Lambda_0}(\chi_{\Lambda_1}^{\q}+\chi_{\Lambda_3}^{\q}+\chi_{\Lambda_4}^{\q})
],\\
s_2&:=\eta^{-8}(-2\chi_{\Lambda_1}+\chi_{\Lambda_3}+\chi_{\Lambda_4}),\\
s_3&:=\eta^{-8}(\chi_{\Lambda_3}-\chi_{\Lambda_4}),\\
s_4&:=\frac{1}{24}\eta^{-16}
[3I^*(d(\eta^{8}{s}_3),d(\eta^{8}{s}_3))
+I^*(d(\eta^{8}{s}_2),d(\eta^{8}{s}_2))].
\end{align}
\end{subequations}
We also put 
$\tilde{s}_{-1}:=s_{-1}$, 
$\tilde{s}_{0}:=s_{0}$, 
$\tilde{s}_{1}:=\eta^{4}s_{1}$, 
$\tilde{s}_{2}:=\eta^{8}s_{2}$, 
$\tilde{s}_{3}:=\eta^{8}s_{3}$, 
$\tilde{s}_{4}:=\eta^{12}s_{4}$.
\begin{rmk}
The functions $s_0,\cdots,s_4$ are reductions of the Jacobi forms 
defined in Section 6 (for the precise statement, see Proposition \ref{631.006} (ii)). 
The functions $\tilde{s}_0,\cdots,\tilde{s}_4$ are twisted by $\eta$ 
so that the automorphic factor do not appear. 
They appear in the calculation of the bilinear forms 
(see Proposition \ref{614.005} (iii)). 
\end{rmk}

\begin{prop}\label{170831.31}
The following results hold. 
\begin{enumerate}
\item $\tilde{s}_0,\cdots,\tilde{s}_4$ give polynomial generators 
of the algebra $S^W$ over $F(\HH) $. 
\item  By $(\mathrm{i})$, we could take $\frac{\p}{\p \tilde{s}_i}$ $(i=-1,\cdots,4)$ 
as $S^W$ free basis of $\Theta_{S^W}$. 
Then $\frac{\p}{\p \tilde{s}_4} \in V$. 
We put $e_1:=\frac{\p}{\p \tilde{s}_4}$. 
\end{enumerate}
\end{prop}
\begin{pf}
We give the proof of (i).  
By Proposition \ref{631.006}, 
$s_0,\cdots,s_4$ are polynomial generators of $S^W$ 
over $F(\HH) $. 
Since $\eta$ is a unit of $F(\HH) $, we have (i). 
For the proof of (ii), by Proposition \ref{631.007}, 
$
\left(\frac{\p}{\p \tilde{s}_4}\right)^2 I^*(d\tilde{s}_i,d\tilde{s}_j)=0
$ 
for $i,j=-1,\cdots,4$. Then by 
Proposition 
\ref{510.003}, we have $\frac{\p}{\p \tilde{s}_4} \in V$. 
\qed
\end{pf}
By Proposition \ref{300.002} and Proposition \ref{170831.31} (ii), 
there exists uniquely the Frobenius structure for $S^W$ 
such that the unit of the product $\circ$ coincides with $e_1$. 
We denote the holomorphic metric for this Frobenius structure by $J_1$. 
Then the dual metric of $J_1$ coincides with $Lie_{e_1}I^*$ 
by Proposition \ref{510.001}. We put $J_1^*:=Lie_{e_1}I^*$. 
\begin{prop}\label{170831.34}
\begin{enumerate}
\item 
For $d\tilde{s}_j$ $(j=-1,\cdots,4)$, we have 
$$\label{520.003}
(J_1^*(d\tilde{s}_i,d\tilde{s}_j))=
\bp 
0&0&0&0&0&1\\
0&6\eta^4 \tilde{E}_4^2&6\eta^4 \tilde{E}_6&0&0&-\frac{5}{6}\eta^4\tilde{E}_4\tilde{s}_1\\
0&6\eta^4 \tilde{E}_6&6\eta^4 \tilde{E}_4&0&0&-\frac{1}{2}\eta^4 \tilde{s}_0\\
0&0&0&12\eta^4&0&0\\
0&0&0&0&4\eta^4&0\\
1&-\frac{5}{6}\eta^4 \tilde{E}_4\tilde{s}_1&-\frac{1}{2}\eta^4\tilde{s}_0&0&0&
\frac{\eta^4}{36}[\frac{2}{3}\tilde{s}_1^2
+\frac{1}{12}\tilde{E}_4\tilde{s}_2^2+\frac{1}{4}\tilde{E}_4\tilde{s}_3^2]
\ep,
$$
where 
\begin{eqnarray*}
E_4&:=&1+240\sum_{n=1}^{\infty}\sigma_3(n)q^n=1+240q+\cdots,\\
E_6&:=&1-504\sum_{n=1}^{\infty}\sigma_5(n)q^n=1-504q-\cdots,
\end{eqnarray*}
$\tilde{E}_{4}:=\eta^{-8}E_4$, 
$\tilde{E}_{6}:=\eta^{-12}E_6$ 
and $\sigma_k(n):=\sum_{d|n}d^k$. 
\item 
The connection matrix of $\nabla$ defined in Theorem \ref{300.002} 
can be calculated 
by $J_1^*(d\tilde{s}_i,d\tilde{s}_j)$ $(i,j=-1,\cdots,4)$ and 
be expressed in terms of $\eta$, 
$\eta'/\eta=\frac{1}{24}E_2$, 
$\tilde{E}_4$, $\tilde{E}_6$ and 
$\tilde{s}_0,\cdots,\tilde{s}_3$. 
%
\end{enumerate}
\end{prop}
\begin{pf}
We give the proof of (i). We first show that $I^*(d\tilde{s}_i,d\tilde{s}_j)$ 
for $i,j=-1,\cdots,4$ could be 
expressed by $\eta,\tilde{E}_i,\tilde{s}_j$ 
$i=2,4,6, j=0,\cdots,4$. 
For the cases $i=-1$ or $j=-1$, then they are shown 
by the equation (\ref{240.002}). 
For the cases $i,j=0,\cdots,4$, then they are shown 
by Proposition \ref{631.008} and the equation (\ref{614.003}). 
Therefore we obtain the results 
because $J_1^*(d\tilde{s}_i,d\tilde{s}_j)=e_1
I^*(d\tilde{s}_i,d\tilde{s}_j)$ for $i.j=-1,\cdots,4$. 
For the proof of (ii), by Theorem \ref{300.002}, the connection $\nabla$ for Frobenius structure 
could be calculated as the Levi-Civita connection of $J_1$, 
we obtain the connection form by the derivatives 
of the $J_1^*(d\tilde{s}_i,d\tilde{s}_j)$ and they could be 
expressed by $\eta,\tilde{E}_i,\tilde{s}_j$ 
$i=2,4,6, j=0,\cdots,4$ using (\ref{641.001}). 
\qed
\end{pf}
\begin{prop}\label{170831.32}
There exist unique 
$v_{0},\cdots,v_{4} \in \Omega_{S^W}$ 
such that 
$v_i-d\tilde{s}_i \in F_0$ for $0 \leq i \leq 3$ and 
$v_4-d\tilde{s}_4 \in F_1$ for the filtration defined in $(\ref{250.001})$
and $\overline{\nabla}v_j=0$ for $0 \leq j \leq 4$, 
where $\overline{\nabla}$ is a connection 
obtained by the reduction of the entries of the connection matrix 
of $\nabla$ defined in Theorem $\ref{300.002}$ 
by $\Omega_{S^W} \to \Omega_{S^W/F(\HH) }$ 
$($more precicely, 
$\overline{\nabla}: \Omega_{S^W} \to \Omega_{S^W/F(\HH) }\otimes_{S^W}
\Omega_{S^W}$ is a composition of the connection 
$\nabla:\Omega_{S^W} \to \Omega_{S^W}\otimes_{S^W} \Omega_{S^W}$ 
with $\Omega_{S^W}\otimes_{S^W} \Omega_{S^W} \to \Omega_{S^W/F(\HH) }\otimes_{S^W} \Omega_{S^W}$ 
induced by the natural morphism 
$\Omega_{S^W} \to \Omega_{S^W/F(\HH) }$$)$. 
We have 
\begin{eqnarray*}\label{520.004}
v_0&=&d\tilde{s}_0-\left((\eta^{2})'\eta^{-2}\tilde{s}_0
-\frac{1}{3}(E_4\eta^{-8})\eta^4\tilde{s}_1\right)\frac{dq}{q},
\\
v_1&=&d\tilde{s}_1-\left((\eta^{2})'\eta^{-2}\tilde{s}_1
-\frac{1}{6}\eta^4\tilde{s}_0\right)\frac{dq}{q},\\
v_2&=&d\tilde{s}_2-\left((\eta^{2})'\eta^{-2}\tilde{s}_2\right)
\frac{dq}{q},\\
v_3&=&d\tilde{s}_3-\left((\eta^{2})'\eta^{-2}\tilde{s}_3\right)
\frac{dq}{q},\\
v_4&=&du_4,
\end{eqnarray*}
where 
\begin{equation}\label{170830.1}
u_4:=\tilde{s}_4+
\frac{-1}{2^9\cdot 3^5}A+\eta^{-4}E_2
\left(\frac{B}{2^9\cdot 3^5}+
\frac{1}{2^4\cdot 3^2}(\tilde{s}_2)^2+\frac{1}{2^4\cdot 3}(\tilde{s}_3)^2\right)
\end{equation}
and 
\begin{eqnarray*}
A&:=&\tilde{E}_6 \tilde{s}_0^2-2\tilde{E}_4^2\tilde{s}_0\tilde{s}_1
+\tilde{E}_4\tilde{E}_6\tilde{s}_1^2,\\
B&:=&
\tilde{E}_4\tilde{s}_0^2-2\tilde{E}_6\tilde{s}_0\tilde{s}_1
+\tilde{E}_4^2\tilde{s}_1^2.
\end{eqnarray*}
Put $v_{-1}:=\frac{1}{2}\frac{dq}{q}$. 
For the $S^W$-free basis $v_{-1},\cdots,v_{4}$, we have 
\begin{equation}\label{520.005}
(\nabla v_{-1},\cdots,\nabla v_{4})
=(v_{-1},\cdots,v_4)
\bp 
0&0&0&0&0&0\\
0&2(\frac{\eta'}{\eta})&-\eta^4/6&0&0&0\\
0&-1/3\eta^4(\tilde{E}_4)&2(\frac{\eta'}{\eta})&0&0&0\\
0&0&0&2(\frac{\eta'}{\eta})&0&0\\
0&0&0&0&2(\frac{\eta'}{\eta})&0\\
0&0&0&0&0&0
\ep
\frac{dq}{q},
\end{equation}
and 
\begin{equation}
\label{520.006}
(J_1^*(v_i,v_j))=
\bp 
0&0&0&0&0&1\\
0&6\eta^4\tilde{E}_4^2&6\eta^4\tilde{E}_6&0&0&0\\
0&6\eta^4\tilde{E}_6&6\eta^4\tilde{E}_4&0&0&0\\
0&0&0&12\eta^4&0&0\\
0&0&0&0&4\eta^4&0\\
1&0&0&0&0&0
\ep.
\end{equation}
\end{prop}
\begin{pf}
By the homogeneity of $J_1^*$, the Levi-Civita connection $\nabla$ 
is homogeneous of degree $0$. 
Thus $\nabla$ preserves the filtration $F_0 \subset F_1 \subset F_2$ 
defined in (\ref{250.001}). 
Then by using Lemma \ref{710.001} inductively, 
we obtain $v_i$ explicitly. 
\qed
\end{pf}

\begin{prop}\label{170831.33}
\begin{enumerate}
\item Put 
\begin{equation}
\label{520.008}
A_0:=
\bp
\eta^{-2}&0\\
0&6\eta^{-2}\\
\ep,\quad
A_1:=
\bp 
\chi_{\Lambda_0}^{\q}&\chi_{\Lambda_1}^{\q}\\
\eta^{-4}(\chi_{\Lambda_0}^{\q})'&\eta^{-4}(\chi_{\Lambda_1}^{\q})'\\
\ep. 
\end{equation}
Then $A_2:=A_0A_1$ satisfies 
\begin{equation}
\label{520.007}
dA_2+\gamma A_2=0, 
\hbox{ for }
\gamma=
\bp 
2(\frac{\eta'}{\eta})&-\eta^4/6\\
-1/3\eta^4(\tilde{E}_4)&2(\frac{\eta'}{\eta})\\
\ep, 
\end{equation} 
and 
\begin{equation}\label{520.012}
(A_2)^T 
\bp 
6\eta^4\tilde{E}_4^2&6\eta^4\tilde{E}_6\\
6\eta^4\tilde{E}_6&6\eta^4\tilde{E}_4
\ep
A_2=\bp 6^3\cdot 48&0\\
0&6^3\cdot 16\ep.
\end{equation}
Put 
$x_{-1}:=v_{-1}, x_2:=\eta^{-2}v_2,x_3:=\eta^{-2}v_3,x_4:=v_4$, 
$
(x_0,x_1):=(v_0,v_1)A_2
$. 
Then the sections $x_{-1},\cdots,x_4$ are flat and we have 
\begin{equation}
\label{520.009}
(J_1^*(x_i,x_j))
=
\bp 
0&0&0&0&0&1\\
0&6^3\cdot 48&0&0&0&0\\
0&0&6^3\cdot 16&0&0&0\\
0&0&0&12&0&0\\
0&0&0&0&4&0&\\
1&0&0&0&0&0
\ep.
\end{equation}
\item Put $\bp y_{-1}&\cdots&y_4\ep:=\bp x_{-1}&\cdots&x_4\ep M^{-1}$ 
for 
\begin{equation}
M:=
\bp 
1&0&0&0&0&0\\
0&72&0&0&0&0\\
0&0&24&-2&0&0\\
0&0&24&1&1&0\\
0&0&24&1&-1&0\\
0&0&0&0&0&1
\ep. 
\end{equation}
Then 
\begin{equation}
\label{520.011}
\bp y_{-1}&\cdots&y_4\ep
=(d(\pi\sqrt{-1}\tau),d(\eta^{-2}\chi_{\Lambda_{0}}),d(\eta^{-2}\chi_{\Lambda_{1}}),
d(\eta^{-2}\chi_{\Lambda_{3}}),d(\eta^{-2}\chi_{\Lambda_{4}}),v_4)
\end{equation}
and 
\begin{equation}
\label{520.010}
(J_1^*(y_i,y_j))=
\bp 
0&0&0&0&0&1\\
0&2&0&0&0&0\\
0&0&2&0&0&0\\
0&0&0&2&0&0\\
0&0&0&0&2&0\\
1&0&0&0&0&0
\ep. 
\end{equation}
in particular, the invariants 
$\pi\sqrt{-1}\tau,
\eta^{-2}\chi_{\Lambda_{0}},
\eta^{-2}\chi_{\Lambda_{1}},
\eta^{-2}\chi_{\Lambda_{3}},
\eta^{-2}\chi_{\Lambda_{4}},
u_4$ 
constitute the flat generator system and 
$\eta^{-2}\chi_{\Lambda_{0}},
\eta^{-2}\chi_{\Lambda_{1}},
\eta^{-2}\chi_{\Lambda_{3}},
\eta^{-2}\chi_{\Lambda_{4}},
u_4$ are polynomial generators of $S^W$ over $F(\HH) $. 
\item The following function 
\begin{eqnarray*}
\F_1&:=&
\frac{1}{2}(\pi \sqrt{-1}\tau)(u_4)^2+\frac{1}{4}u_4(b_0^2+b_1^2+b_2^2+b_3^2)\\
&&+f_0(b_0b_1b_2b_3)\\
&&+\frac{1}{4}f_1(b_0^4+b_1^4+b_2^4+b_3^4)\\
&&+\frac{1}{6}
f_2(b_0^2b_1^2+b_0^2b_2^2+b_0^2b_3^2+b_1^2b_2^2+b_1^2b_3^2+b_2^2b_3^2)
\end{eqnarray*}
gives a potential for the above Frobenius structure where 
$f_0,f_1,f_2$ are defined in Theorem $\ref{400.001}$. 
\item The function $u_4$ defined $(\ref{170830.1})$ 
coincides with the function $b_4$ in Theorem $\ref{400.001} (\mathrm{i})$.
\end{enumerate}
\end{prop}
\begin{pf}
We give the proof of (i). 
By Proposition \ref{641.002} (i), (iii), we obtain (\ref{520.007}) and (\ref{520.012}) 
respectively. 
By (\ref{520.007}), $x_i$ are flat. 
By (\ref{520.012}), we have (\ref{520.009}). 
For the proof of (ii), by the fact that the determinant of $A_1$ is $4$ 
(Proposition \ref{641.002} (ii)), 
we have (\ref{520.011}).
We give the proof of (iii). 
In the proof of Proposition \ref{170831.34} (i), we obtain $I^*(d\tilde{s}_i,d\tilde{s}_j)$. 
We could express $I^*(y_i,y_j)$ by $\eta,\eta^4\chi_{\Lambda_0}, \eta^4\chi_{\Lambda_1}, 
\eta'/\eta, \tilde{E}_i,y_j$ 
$j=0,\cdots,4$ using the results of Proposition \ref{641.003}. 
Then by Proposition \ref{510.004}, we have (iii). By the result of (iii), 
we obtain $u_4=b_4$, which gives the proof of (iv). 
\qed
\end{pf}

\begin{pf}[of Theorem \ref{400.001}]
By Theorem \ref{300.002}, the connection in Theorem \ref{400.001} 
is the same as one in Proposition \ref{170831.34} (ii). 
Then (i) is a consequence of Proposition \ref{170831.33} (ii) and (iv). 
By Proposition \ref{170831.33} (iv), we have $e_0=e_1$. 
Then we have $J_0=J_1$ and the function $\F_1$ gives the potential 
of the Theorem \ref{400.001}. 
\qed\end{pf}

\section{Jacobi forms}

\subsection{Domains and differential operators}
In order to define Jacobi forms, we prepare domains, group actions 
on them and differential operators. 

We formulate Jacobi forms on the new domains (see (\ref{170831.22}))  
which have 1 extra parameter 
than the domain $Y$ so that the automorphic factors for the group actions do not appear. 
We also formulate differential operators which 
are invariant for the group actions on the new domains (see (\ref{170831.23})). 
We use notations introduced in Section 2. 

\subsubsection{Domains and group actions}\label{170831.21}
We put 
$U:=\{q \in \CC\,|\,|q|<1\}$, 
$Y_{\HH}:=\HH \times \mathfrak{h}_{\CC} \times \CC$, 
$Y_{U}:=U \times \mathfrak{h}_{\CC} \times \CC$, 
$Y_{\{0\}}:=\{0\} \times \mathfrak{h}_{\CC} \times \CC$. 
We prepare $\mathbb{C}^*$ with the coordinate $\omega$. 
The maps: $\HH \to U,\ \tau \mapsto e^{2\pi\sqrt{-1}\tau}$ and 
$\{0\} \to U$ induce the following diagram:
\begin{equation}\label{611.001}
\begin{CD}
\CC^* \times Y_{\HH} 
@>\varphi_1>> 
\CC^*\times Y_{U} 
@<\varphi_2<<  
\CC^*\times Y_{\{0\}}
\\
@V{\pi}VV @V{\pi}VV @V{\pi}VV\\
\CC^* \times \HH 
@>\varphi_1>> 
\CC^*\times U
@<\varphi_2<<  
\CC^*\times\{0\}
\end{CD},
\end{equation}
where $\pi$ denotes natural projection. 

We shall define the Metaplectic group (cf. \cite[p.253]{Kac}). 
We first introduce $SL_2(\ZZ)$ action on $\HH$ 
by $A\cdot \tau:=\frac{a\tau+b}{c\tau+d}$ 
for $A=\bp a&b\\c&d\ep \in SL_2(\ZZ)$ and $\tau \in \HH$. 
We define the Metaplectic group $Mp_2(\ZZ)$ by the set of pair 
$$
\left\{\left(\bp a&b\\c&d\ep,j\right)\in SL_2(\ZZ)\times F(\HH) \,|\,j^2=c\tau+d
\right\}
$$ with the product. 
$(A_1,j_1)(A_2,j_2)=(A_1A_2,j_3)$ where $j_3(\tau)=j_1(A_2\cdot \tau)j_2(\tau)$. 
We define the action of $Mp_2(\ZZ)$ on $\CC^*\times \HH$ by 
$(A,j)\cdot (\omega,\tau)=(j\omega,A\cdot \tau)$. \\
We define the action of $Mp_2(\ZZ)$ on $(\omega,\tau,z,t) \in \CC^*\times Y_{\HH}$ 
by 
$$
(A,j)\cdot (\omega,\tau,z,t)=
\left(j\omega,A\cdot \tau, \frac{z}{c\tau+d},
t+\frac{c\langle z,z\rangle}{2(c\tau+d)}\right) \in \CC^* \times Y_{\HH}
$$
for $A=\bp a&b\\c&d\ep \in SL_2(\ZZ)$ 
and $\langle z,z\rangle$ is a non-degenerate 
$\CC$-bilinear form on $h$ induced by a non-degenerate 
$\RR$-bilinear form on $F_{\mathrm{cl}}$ induced by $(F,I)$. 

We define the $W$ action on $Y_{\HH}$ 
by the isomorphism $Y \simeq Y_{\HH}$ 
which was obtained in (\ref{240.003}) and $W$ action on $Y$ 
defined in Section \ref{240.001}. 
The domain $\CC^* \times Y_{\HH}$ also has the $W$ action 
whose action on $\CC^*$ is trivial. 

\subsubsection{Differential operators}
We remind the notion of ``$f$-related" which we shall use frequently.  
Let $f:M_1 \to M_2$ be a morphism of complex manifolds $M_1$ and $M_2$. 
Let $X_i$ be a holomorphic vector field on $M_i$ ($i=1,2$). 
The vector fields $X_1$ and $X_2$ are called $f$-related if 
$df((X_1)_p)=(X_2)_{f(p)}$ for any $p \in M_1$ 
where $df:T_pM_1 \to T_{f(p)}M_2$. 
Then we have a commutative diagram:
\begin{equation}
\begin{CD}
F(M_1)  @<{f^*}<< F(M_2)\\
@V{X_1}VV @V{X_2}VV \\
F(M_1)  @<{f^*}<< F(M_2)
\end{CD}.
\end{equation}
for the function spaces $F(M_i)$ ($i=1,2$). 
For this case, we use the same notation $X_1=X_2=X$. 
If $P_i$ is a differential operator on $M_i$ ($i=1,2$) such that 
the diagram 
\begin{equation}
\begin{CD}
F(M_1)  @<{f^*}<< F(M_2)\\
@V{P_1}VV @V{P_2}VV \\
F(M_1)  @<{f^*}<< F(M_2)
\end{CD}
\end{equation}
commutes, then we also call $D_1$ and $D_2$ ``$f$-related" and 
use the same notation $P=P_1=P_2$. 

We return to the situation of Section \ref{170831.21}. 
On each space of the diagram (\ref{611.001}), we define the vector field 
$\omega\frac{\p}{\p \omega}$. They are $\varphi_i$-related ($i=1,2$) and $\pi$-related, 
thus we simply denote them by $E_{\omega}$. 

On each space of the upper line of the diagram (\ref{611.001}), 
we define the vector field $\frac{\p}{\p \Lambda_0}$. 
They are $\varphi_i$-related ($i=1,2$), thus we simply denote them by $E$. 

On each space 
$\CC^* \times Y_{\HH}$, 
$\CC^* \times Y_{U}$ and 
$\CC^* \times Y_{\{0\}}$, 
we define the differential operators  
\begin{eqnarray*}
&&\omega^{-4}
\left[D+2\frac{\eta'}{\eta}
\{l\frac{\p}{\p \Lambda_0}+
\omega \frac{\p^2}{\p \Lambda \p \omega}\}
\right],\\
&&\omega^{-4}
\left[2q\frac{\p^2}{\p q \p \Lambda_0}-\overline{D}+\frac{2}{24}E_2
\{l\frac{\p}{\p \Lambda_0}+
\omega \frac{\p^2}{\p \Lambda \p \omega}\}
\right],\\
&&\omega^{-4}
\left[-\overline{D}+\frac{2}{24}\{l\frac{\p}{\p \Lambda_0}+
\omega \frac{\p^2}{\p \Lambda \p \omega}\}
\right]
\end{eqnarray*}
respectively, 
where the symbol $'$ is defined in Section \ref{170930}. 
They are $\varphi_i$-related ($i=1,2$), thus we simply denote them by $\mathcal{D}$. 

On each space $\CC^*\times \HH$, $\CC^*\times U$ and $\CC^*\times \{0\}$, 
we define vector fields 
\begin{eqnarray*}
&&\omega^{-4}
\left[
\frac{1}{24}E_2(e^{2\pi\sqrt{-1}\tau})
\omega\frac{\p}{\p \omega}+\frac{1}{2\pi\sqrt{-1}}\frac{\p}{\p \tau}
\right],\\
&&\omega^{-4}[\frac{1}{24}E_2(q)
\omega\frac{\p}{\p \omega}
+q\frac{\p}{\p q}],\\
&&\omega^{-4}[\frac{1}{24}\omega\frac{d}{d \omega}]
\end{eqnarray*}
respectively. 
They are $\varphi_1$-related and $\varphi_2$-related, thus 
we simply denote them by $\delta_q$. 
\begin{prop}
\begin{enumerate}
\item
On the space $\CC^*\times Y_{\HH}$, 
$E_{\omega}, E$ and $\mathcal{D}$ are invariant w.r.t. 
the $Mp_2(\ZZ)$ action and $W$ action. 
We also have $[E_{\omega},\mathcal{D}]=(-4)\mathcal{D}$. 
\item
On the space $\CC^*\times \HH$, 
$E_{\omega}$ and $\delta_q$ are invariant w.r.t. the $Mp_2(\ZZ)$ action. 
We also have $[E_{\omega},\delta_q]=(-4)\delta_q$. 
\end{enumerate}
\end{prop}
We could easily check the above facts so we omit a proof.

\subsubsection{Function spaces}

For $k \in \frac{1}{2}\ZZ$, $m \in \ZZ$ and for each $X=\HH$, $U$ and $\{0\}$, we put 
\begin{eqnarray*}
F_{k,m}(\CC^* \times Y_X)
&:=&
\{f:\CC^* \times Y_X \to \CC:holomorphic \,|\,
E_{\omega}f=(-2k)f,\ 
Ef=mf\},\\
F_k(\CC^*\times X)
&:=&\{f:\CC^*\times X \to \CC:holomorphic \,|
E_{\omega}f=(-2k)f\}.
\end{eqnarray*}
Since the vector fields $E_{\omega}$ and $E$ 
are $\varphi_i$-related ($i=1,2$), we have the induced morphisms $\varphi_i^*$ 
on these function spaces. 

By $[E_{\omega},\mathcal{D}]=(-4)\mathcal{D}$ and 
$[E_{\omega},\delta_q]=(-4)\delta_q$, the differential operators 
$\mathcal{D}$, $\delta_q$ define 
\begin{eqnarray}
\mathcal{D}&:&F_{k,m}(\CC^* \times Y_X)
\to 
F_{k+2,m}(\CC^* \times Y_X),\\
\label{170831.23}
\delta_q&:&F_{k}(\CC^*\times X) \to F_{k+2}(\CC^*\times X).
\end{eqnarray}
Since the differential operator $\mathcal{D}$ and $\delta_q$ are $\varphi_i$-related ($i=1,2$), 
the above action on the function spaces commutes 
with $(\varphi_1)^*$ and $(\varphi_2)^*$ 
(cf. Proposition \ref{614.005}). 

\begin{prop}
\begin{enumerate}
\item
For the natural inclusion mapping:
$
\pi^*:
F_k(\CC^* \times X) 
\to 
F_{k,0}(\CC^* \times Y_X)
$, 
we easily check that 
\begin{equation}\label{613.001}
\mathcal{D}({f}{F})
=2\delta_q({f})E({F})
+{f}\mathcal{D}({F})
\end{equation}
for ${f} \in F_k(\CC^* \times X) $, 
${F} \in F_{k',m}(\CC^* \times Y_X)$. 
\item
For $\eta/\omega \in F_1(\CC^*\times \HH)$, 
\begin{equation}\label{613.003}
\delta_q(\eta/\omega)=0.
\end{equation}
\end{enumerate}
\end{prop}
We could easily check the above facts so we omit a proof.

For $F_1 \in F_{k,m}(\CC^* \times Y_X)$, 
$F_2 \in F_{k',m'}(\CC^* \times Y_X)$, we put 
\begin{equation}\label{613.002}
\mathcal{I}(F_1,F_2)
:=
\frac{1}{2}
\left[
\mathcal{D}
(F_1F_2)
-
\mathcal{D}
(F_1)F_2
-
\mathcal{D}
(F_2)F_1
\right]\in 
F_{k+k'+2,m+m'}(\CC^* \times Y_X).
\end{equation}
We remark that 
$\varphi_1^*$ is injective 
and 
$\varphi_2^*$ is surjective,  
because $\varphi_1$ is dominant 
and $\varphi_2$ is a closed immersion respectively. 
Hereafter we frequently identify 
$F_{k,m}(\CC \times Y_U)$ with the subspace 
$\varphi_1^*(F_{k,m}(\CC \times Y_U))$. 
For $f \in F_{k,m}(\CC \times Y_U)$, 
we call $\varphi_2^*(f)$ the initial term of $f$ 
(see (\ref{631.002})). 
\subsubsection{Correspondence of the function spaces}
\label{614.004}
By the natural projection $\CC^* \times Y_{\HH} \to Y_{\HH} \simeq Y$ 
(resp. $\CC^*\times \HH \to \HH$), 
we have $F_m(Y) \to F_{0,m}(\CC^*\times Y_{\HH})$ 
(resp. $F(\HH)  \to F_0(\CC^*\times \HH)$). 
We also define their twisted version. 
\begin{defn}
For each $k \in \frac{1}{2}\ZZ$, we define the isomorphism 
\begin{eqnarray*}
L_k&:&F_m(Y) \to F_{k,m}(\CC^*\times Y_{\HH}),\\
L_k&:&F(\HH)  \to F_k(\CC^*\times \HH)
\end{eqnarray*}
by $f \mapsto \omega^{-2k}f$. We put $\hat{E}_k:=L_{k}E_k$ for $k=2,4,6$. 
\end{defn}

\begin{prop} \label{614.005}
\begin{enumerate}
\item We have 
\begin{equation}\label{614.001}
\begin{CD}
F(\HH)  @>{L_k}>> F_k(\CC^*\times \HH)
@<{\varphi_1^*}<< F_k(\CC^* \times U) @>{\varphi_2^*}>> F_k(\CC^*\times \{0\})\\
@V{\p_k}VV @V{\delta_q}VV @V{\delta_q}VV @V{\delta_q}VV 
\\
F(\HH)  @>{L_{k+2}}>> F_{k+2}(\CC^*\times \HH)
@<{\varphi_1^*}<< F_{k+2}(\CC^* \times U) @>{\varphi_2^*}>> F_{k+2}(\CC^*\times \{0\})\\
\end{CD},
\end{equation}
where $\p_k:F(\HH)  \to F(\HH) $ by 
$\p_k(f)=f'-\frac{k}{12}E_2 f=\eta^{2k}(\eta^{-2k}f)'$
 for $f \in F(\HH) $. 
\item We have 
\begin{equation}\label{614.002}
\begin{CD}
F_m(Y) 
@>{L_{k}}>> 
F_{k,m}(\CC^*\times Y_{\HH})
@<{\varphi_1^*}<< 
F_{k,m}(\CC^*\times Y_U) 
@>{\varphi_2^*}>> 
F_{k,m}(\CC^*\times Y_{\{0\}})\\
@V{\mathcal{D}_k}VV @VV{\mathcal{D}}V @VV{\mathcal{D}}V @VV{\mathcal{D}}V 
\\
F_m(Y) @>{L_{k+2}}>> 
F_{k+2,m}(\CC^*\times Y_{\HH})
@<{\varphi_1^*}<< 
F_{k+2,m}(\CC^*\times Y_U) 
@>{\varphi_2^*}>> 
F_{k+2,m}(\CC^*\times Y_{\{0\}})\\
\end{CD},
\end{equation}
where 
$\mathcal{D}_k(F):=[2(l-2k)\frac{\eta'}{\eta}\frac{\p}{\p t}+D](F)
=\frac{1}{\eta^{l-2k}}D(\eta^{l-2k}F)$ 
for $F \in F_m(Y)$. 
\item For $F_1 \in F_{m}(Y)$, 
$F_2 \in F_{m'}(Y)$, 
\begin{equation}\label{614.003}
I^*(d(\eta^{-2k}F_1),d(\eta^{-2k'}F_2))
=
\eta^4\left(\frac{\eta}{\omega}\right)^{-2k-2k'-4}
\mathcal{I}(L_{k}(F_1),L_{k'}(F_2)).
\end{equation}
\end{enumerate}
\end{prop}
\begin{pf}
We have (i) and (ii) by direct calculations. 
For (iii), we have 
\begin{eqnarray*}
&&\frac{1}{\omega^{-2k-2k'-4}}
\mathcal{I}
(\omega^{-2k}F_1,\omega^{-2k'}F_2)\\
&=&
\frac{1}{2}
\frac{1}{\omega^{-2k-2k'-4}}
\left[
\mathcal{D}(\omega^{-2k-2k'}F_1F_2)
-\mathcal{D}(\omega^{-2k}F_1)\omega^{-2k'}F_2
-\mathcal{D}(\omega^{-2k'}F_2)\omega^{-2k}F_1
\right]\\
&=&
\frac{1}{2}
\left[
\frac{1}{\eta^{l-2k-2k'}}D(\eta^{l-2k-2k'}F_1F_2)
-\frac{1}{\eta^{l-2k}}D(\eta^{l-2k}F_1)F_2
-\frac{1}{\eta^{l-2k'}}D(\eta^{l-2k'}F_2)F_1
\right]\\
&=&
\frac{1}{2}
\frac{1}{\eta^{l-2k-2k'}}
\left[
D(\eta^{l-2k-2k'}F_1F_2)
-D(\eta^{l-2k}F_1)(\eta^{-2k'}F_2)
-D(\eta^{l-2k'}F_2)(\eta^{-2k}F_1)
\right]\\
&=&
\frac{1}{2}
\frac{1}{\eta^{-2k-2k'}}
\left[
D(\eta^{-2k-2k'}F_1F_2)
-D(\eta^{-2k}F_1)(\eta^{-2k'}F_2)
-D(\eta^{-2k'}F_2)(\eta^{-2k}F_1)
\right]\\
&=&
\frac{1}{\eta^{-2k-2k'}}
I^*(d(\eta^{-2k}F_1),d(\eta^{-2k'}F_2)).
\end{eqnarray*}
Here we used (ii) and (\ref{250.002}) and (\ref{250.003}).
\qed\end{pf}

\subsection{Definitions of modular forms and Jacobi forms}
\subsubsection{Spaces of modular forms and Jacobi forms}

We put 
\begin{eqnarray}\label{621.001}
&&F_{k,m}^W(\CC^*\times Y_{\HH}):=\{
f \in F_{k,m}(\CC^*\times Y)
\,|\,
f \hbox{ is invariant by the action of }W\}, \nonumber\\
&&F_{k,m}^W(\CC^*\times Y_U):=\{
f \in F_{k,m}(\CC^*\times Y_U)
\,|\,
\varphi_1^*(f) \in F_{k,m}^W(\CC^*\times Y_{\HH})
\}, \nonumber\\
&&J_{k,m}:=\{
f \in F_{k,m}^W(\CC^*\times Y_U)
\,|\,
\varphi_1^*(f) 
\hbox{ is invariant by the action of }
Mp_2(\ZZ)
\}, \\
\label{170831.22}
&&M_k:=\{f \in F_k(\CC^*\times U)\,|\,\varphi_1^*(f) 
\hbox{ is invariant by the action of }
Mp_2(\ZZ)\}.
\end{eqnarray}
We call $M_k$ and $J_{k,m}$ 
the space of the modular form of weight $k$ and 
the space of Jacobi forms of weight $k$ and index $m$ 
respectively. 

We put $M_k^{\omega=1}:=(L_k)^{-1}(\varphi_1^*(M_{k}))$ and 
$J_{k,m}^{\omega=1}:=(L_{k})^{-1}(\varphi_1^*(J_{k,m}))$. 
We remark that $M_k^{\omega=1}$ is a usual space of 
modular forms of weight $k$ on $\HH$. 
Also we remark that the space $J_{k,m}^{\omega=1}$ for $A_1^{(1,1)}$ case 
corresponds to the space of the even weak Jacobi forms defined in Eichler-Zagier \cite{EZ} 
by restricting $\HH \times \mathfrak{h}_{\CC} \times \CC$ 
to $\HH \times \mathfrak{h}_{\CC}\times \{0\}$. 


\begin{prop}
\begin{enumerate}
\item We have 
\begin{equation}\label{621.002}
\begin{CD}
M_k^{\omega=1} @>{\sim}>{L_k}> \varphi_1^*(M_k) 
@<{\sim}<{\varphi_1^*}< M_k 
@>>{\varphi_2^*}> M_{k}^0\\
@V{\p_{k}}VV @V{\delta_q}VV @V{\delta_q}VV @V{\delta_q}VV\\
M_{k+2}^{\omega=1} @>{\sim}>{L_{k+2}}> \varphi_1^*(M_{k+2})
@<{\sim}<{\varphi_1^*}< M_{k+2}
 @>>{\varphi_2^*}> M_{k+2}^0
\end{CD}, 
\end{equation}
where we put $M_k^0:=\varphi_2^*(M_k)$ and call it 
the space of the initial terms of the modular forms of weight $k$. 
\item We have 
\begin{equation}\label{621.003}
\begin{CD}
J_{k,m}^{\omega=1} @>{\sim}>{L_k}> \varphi_1^*(J_{k,m})
@<{\sim}<{\varphi_1^*}< J_{k,m} 
@>>{\varphi_2^*}> J_{k,m}^0\\
@V{\mathcal{D}_k}VV @V{\mathcal{D}}VV @V{\mathcal{D}}VV @V{\mathcal{D}}VV \\
J_{k+2,m}^{\omega=1} @>{\sim}>{L_{k+2}}> \varphi_1^*(J_{k+2,m})
@<{\sim}<{\varphi_1^*}< J_{k+2,m} 
@>>{\varphi_2^*}> J_{k+2,m}^0\\
\end{CD}, 
\end{equation}
where we put 
$J_{k,m}^0:=\varphi_2^*(J_{k,m})$ and call it 
the space of the initial terms of the Jacobi forms 
of weight $k$ and index $m$. 
\end{enumerate}
\end{prop}
\begin{pf}
Since $\delta_q$ 
is $Mp_2(\ZZ)$-equivariant, the diagram (\ref{614.001}) 
gives the commutative diagram (\ref{621.002}). 
Also since $\mathcal{D}$ is $Mp_2(\ZZ)$-equivariant, 
the diagram (\ref{614.002}) gives the commutative diagram (\ref{621.003}).  
\qed
\end{pf}

\subsubsection{Determinant construction of Jacobi forms}
\label{622.002}
\begin{prop}\label{622.001}
Let 
$f_1,\cdots,f_n \in F_{k,m}^W(\CC^* \times Y_{\HH})$, 
$\overline{f_1},\cdots,\overline{f_n} \in 
F_{k,m}^W(\CC^*\times Y_U)$, 
$F_1,\cdots,F_n \in F_{k}(\CC^*\times \HH)$, 
$\overline{F_1},\cdots,\overline{F_n} \in F_{k}(\CC^*\times U)$, 
$\gamma_1,\cdots,\gamma_n \in \QQ$ 
and $\rho:Mp_2(\mathbb{Z}) \to GL_n(\mathbb{C})$ 
satisfy the following conditions:
\begin{enumerate}
\item 
$f_1=q^{\gamma_1}\overline{f_1}(\omega,q,z,t), 
\cdots,
f_n=q^{\gamma_n}\overline{f_n}(\omega,q,z,t), 
$
\item 
$F_1(\omega,\tau)=q^{\gamma_1}\overline{F_1}(\omega,q), 
\cdots,
F_n(\omega,\tau)=q^{\gamma_n}\overline{F_n}(\omega,q), 
$
\item  
$$
(f_1(g\cdot(\omega,\tau,z,t)),\cdots,
f_n(g\cdot(\omega,\tau,z,t)))
=
(f_1(\omega,\tau,z,t),\cdots,
f_n(\omega,\tau,z,t))
\rho(g),
$$
for $g \in Mp_2(\ZZ)$, 
\item 
$$
(F_1(g\cdot(\omega,\tau)),\cdots,
F_n(g\cdot(\omega,\tau)))
=
(F_1(\omega,\tau),\cdots,
F_n(\omega,\tau))
\rho(g). 
$$
\end{enumerate}
Then we have $24\gamma \in \mathbb{Z}$ 
for $\gamma=\gamma_1+\cdots+\gamma_n$ and for 
$$
\Delta:=
\bp 
f_1&\cdots&f_n\\
\delta_q^{\beta_2}F_1&\cdots&\delta_q^{\beta_2}F_n\\
\delta_q^{\beta_3}F_1&\cdots&\delta_q^{\beta_3}F_n\\
\cdots&\cdots&\cdots\\
\delta_q^{\beta_n}F_1&\cdots&\delta_q^{\beta_n}F_n
\ep, 
$$
$(\eta/\omega)^{-24\gamma}\det \Delta$ is a Jacobi form of type $X_l^{(1,1)}$ 
of weight $\beta-12\gamma$ and index $m$ 
for $\beta=kn+2\{\beta_2+\cdots+\beta_n\}$.  
\end{prop}
\begin{pf}
We first study the behavior of $\det \Delta$ at $q=0$. 
We first observe that 
for $i \geq 2$, the $(i,j)$ component of $\Delta$ 
is $\delta_q^{\beta_i}(q^{\gamma_j}\overline{F_j})$ 
, so it could be factored as 
$q^{\gamma_j}b_{i,j}(q)$ 
by some holomorphic function $b_{i,j}(q)\in F_{k+2\beta_i}(\CC^*\times U)$. 
Thus $\det \Delta$ is factored 
as 
$\det \Delta=q^{\gamma}F(\omega,q,z,t)$ 
by some holomorphic function $F(\omega,q,z,t)\in 
F_{\beta,m}^W(\CC^*\times Y_U)$. 

Second, we study the modular invariance 
of $\det \Delta$. 
Here we consider the group of characters:
$G:=\{\chi:Mp_2(\ZZ) \to \CC^*\,|\,\chi \hbox{ is a group homomorphism}\}$. 
We call $f \in F_{k,m}^W(\CC^* \times Y_{\HH})$ quasi-invariant (resp. invariant) 
for $Mp_2(\ZZ)$ if 
there exist $\chi \in G$ s.t. 
$f(g\cdot (\omega,\tau,z,t))=\chi(g)f(\omega,\tau,z,t)$ 
(resp. $f(g\cdot (\omega,\tau,z,t))=f(\omega,\tau,z,t)$) 
for all $g \in Mp_2(\ZZ)$. 

We show that $\det \Delta$ is quasi-invariant. 
Since $\delta_q$ is invariant w.r.t. $Mp_2(\mathbb{Z})$, 
we have 
$$
((\delta_q^k F_1)(g\cdot(\omega,\tau)),\cdots,
(\delta_q^k F_n)(g\cdot(\omega,\tau)))
=
(\delta_q^k F_1(\omega,\tau),\cdots,
\delta_q^k F_n(\omega,\tau))
\rho(g),
$$
thus we have 
$$
\Delta(g\cdot(\omega,\tau,z,t))
=
\Delta(\omega,\tau,z,t) \rho(g). 
$$
Then we obtain the following relation 
of the determinants of both sides:
$$
\det (\Delta(g\cdot(\omega,\tau,z,t)))
=\det \rho(g)
\det (\Delta(\omega,\tau,z,t)).
$$
Thus $\det \Delta$ is quasi-invariant for $Mp_2(\ZZ)$. 

We show $24\gamma \in \ZZ$. 
We take an element $T:=(\bp 1&1\\0&1\ep,1) \in Mp_2(\ZZ)$. 
Then by $\det \Delta=q^{\gamma}F(\omega,q,z,t)$, we have 
$\det(\rho(T))=e^{2\pi\sqrt{-1}\gamma}$. 
We remind the fact that the character group $G$ is isomorphic 
to cyclic group of order $24$. 
Then $(\det(\rho(T)))^{24}=(e^{2\pi\sqrt{-1}\gamma})^{24}$ 
must be $1$. Thus we have $24\gamma \in \ZZ$. 

Put $\varphi:=(\eta/\omega)^{-24\gamma}\det \Delta$. 
Since $\eta/\omega$ is quasi-invariant for $Mp_2(\ZZ)$, 
$\varphi$ is also quasi-invariant for $Mp_2(\ZZ)$. 
We see that $\varphi$ is invariant w.r.t. $T \in Mp_2(\ZZ)$. 
Since the elements of the character group $G$ is distinguished 
by the value on $T=(\bp 1&1\\0&1\ep,1) \in Mp_2(\ZZ)$, 
we see that $\varphi$ is invariant for $Mp_2(\ZZ)$. 

Then $(\eta/\omega)^{-24\gamma} \det \Delta$ has a modular property 
and satisfy cusp conditions, so it satisfies 
the conditions of Jacobi forms. 
\qed\end{pf}

\subsubsection{Theorem of Wirthm\"uller}
We put $J_{*,*}:=\bigoplus_{k,m \in \ZZ}J_{k,m}$, $M_*:=\bigoplus_{k \in \ZZ}M_k$. 
We review the result of Wirthm\"uller \cite{Chevalley5}. 
\begin{thm}[(3.6) Theorem in \cite{Chevalley5}]\label{623.001}
\begin{enumerate}
\item For the Jacobi forms of types 
$X_l^{(1,1)} (X_l=A_l,B_l,\cdots,G_2)$ except $E_8^{(1,1)}$, 
there exist Jacobi forms $\phi_i \in J_{k_i,m_i}(i=0,\cdots,l)$ 
for some $k_i,m_i \in \ZZ$ such that 
these are polynomial generator of $J_{**}$ over $M_{*}$. 
\item The functions $L^{-1}_{k_0}(\varphi_1^*(\phi_0)),\cdots,
L^{-1}_{k_l}(\varphi_1^*(\phi_l))$ are polynomial generator 
of $S^W$ over $F(\HH) $. 
\end{enumerate}
\end{thm}
\begin{rmk}
For the Jacobi forms of type $E_8^{(1,1)}$, 
the existence of the polynomial generator is not known. 
\end{rmk}

\subsection{$D_4^{(1,1)}$ case}

Hereafter we restrict ourselves to the case 
of type $D_4^{(1,1)}$. 
Then Theorem \ref{623.001} for this case 
asserts the existence of the Jacobi forms 
$\phi_0 \in J_{0,1}$, 
$\phi_1 \in J_{-2,1}$, 
$\phi_2 \in J_{-4,1}$, 
$\phi_3 \in J_{-4,1}$, 
$\phi_4 \in J_{-6,2}$ 
such that these are polynomial generator of 
$J_{*,*}$ over $M_{*}$. 
Wirthm\"uller constructs them by the technique of lifting from lower rank Jacobi forms. 

Here we shall construct Jacobi forms 
by the method of ``the determinant construction" introduced 
in Section \ref{622.002} 
using the normalized characters introduced in 
(\ref{260.002}). 

For $i=0,1,3,4$, put $\hat{\chi}_{\Lambda_i}:=L_0(\chi_{\Lambda_i}) 
\in F_{0,1}(\CC^* \times Y_{\HH})$, 
$\hat{ch}_{\Lambda_i}:=L_0(ch_{\Lambda_i})
\in F_{0,1}(\CC^* \times Y_{\HH})$,  
$\hat{\chi}_{\Lambda_i}^{\q}:=L_0(\chi_{\Lambda_i}^{\q})
\in F_{0}(\CC^* \times {\HH})$, 
$\hat{ch}_{\Lambda_i}^{\q}:=L_0(ch_{\Lambda_i}^{\q})
\in F_{0}(\CC^* \times {\HH})$. 

\begin{prop}\label{631.001}
\begin{enumerate}
\item For $i=0,1,3,4$, we have 
\begin{eqnarray*}\label{631.002}
\hat{ch}_{\Lambda_i}&\in& F_{0,1}(\CC^* \times Y_{U}),  \\
\hat{ch}_{\Lambda_i}^{\q}&\in& F_{0}(\CC^* \times {U}), \\
\hat{\chi}_{\Lambda_i}&=&q^{m_{\Lambda_i}}\hat{ch}_{\Lambda_i},\\
\hat{\chi}_{\Lambda_i}^{\q}&=&q^{m_{\Lambda_i}}\hat{ch}_{\Lambda_i}^{\q},
\end{eqnarray*}
where $m_{\Lambda_i}$ are defined in Section $\ref{260.002}$ 
and for our case we have $m_{\Lambda_0}=-\frac{4}{24}$, 
$m_{\Lambda_1}=m_{\Lambda_3}=m_{\Lambda_4}=\frac{8}{24}$. 
\item 
\begin{align}
\bp 
\hat{\chi} _{\Lambda_0}^T\\
\hat{\chi} _{\Lambda_1}^T\\
\hat{\chi} _{\Lambda_3}^T\\
\hat{\chi} _{\Lambda_4}^T
\ep
=
e^{2\pi\sqrt{-1}/3}
\bp
-1&0&0&0\\
0&1&0&0\\
0&0&1&0\\
0&0&0&1
\ep
\bp 
\hat{\chi} _{\Lambda_0}\\
\hat{\chi} _{\Lambda_1}\\
\hat{\chi} _{\Lambda_3}\\
\hat{\chi} _{\Lambda_4}
\ep,\, 
\bp 
\hat{\chi} _{\Lambda_0}^S\\
\hat{\chi} _{\Lambda_1}^S\\
\hat{\chi} _{\Lambda_3}^S\\
\hat{\chi} _{\Lambda_4}^S
\ep
=\frac{1}{2}\bp
1&1&1&1\\
1&1&-1&-1\\
1&-1&1&-1\\
1&-1&-1&1
\ep
\bp 
\hat{\chi} _{\Lambda_0}\\
\hat{\chi} _{\Lambda_1}\\
\hat{\chi} _{\Lambda_3}\\
\hat{\chi} _{\Lambda_4}
\ep,
\end{align}
where $\hat{\chi} _{\Lambda_i}^g:=\hat{\chi} _{\Lambda_i}(g\cdot (\omega,\tau,z,t))$ 
for $g \in Mp_2(\ZZ)$ $(i=0,1,3,4)$ 
and 
$
T=(\bp 1&1\\0&1\ep,1),\ 
S=(\bp 0&-1\\1&0\ep,\sqrt{\tau}) 
\in Mp_2(\ZZ)$ 
and $\sqrt{\tau}$ is defined as 
$0 <\arg (\sqrt{\tau})<\frac{\pi}{2}$. 
\item For $i=0,1,3,4$, the initial terms of $\hat{ch}_{\Lambda_i}$ are 
\begin{eqnarray*}\label{631.002}
\varphi_2^*(\hat{ch}_{\Lambda_0})&=&e^{\Lambda_0}S(0),\\ 
\varphi_2^*(\hat{ch}_{\Lambda_1})&=&e^{\Lambda_0}S(\omega_1), \\
\varphi_2^*(\hat{ch}_{\Lambda_3})&=&e^{\Lambda_0}S(\omega_3), \\
\varphi_2^*(\hat{ch}_{\Lambda_4})&=&e^{\Lambda_0}S(\omega_4), 
\end{eqnarray*}
where $\omega_k$ is a fundamental weight, 
\begin{equation}\label{631.003}
S(\omega_k):=\sum_{w \in W_{\mathrm{cl}}/\hbox{the isotropy subgroup of }\omega_k}
e^{w\cdot \omega_k}, 
\end{equation}
where $W_{\mathrm{cl}}$ is a finite Weyl group generated by 
$w_{\alpha_1},\cdots,w_{\alpha_4}$. 
\item $\hat{\chi} _{\Lambda_1}^{\q}=\hat{\chi} _{\Lambda_3}^{\q}=\hat{\chi} _{\Lambda_4}^{\q}$. 
\item $\mathcal{D}(\hat{\chi} _{\Lambda_i})=0$ for $i=0,1,3,4$. 
\item For $\displaystyle{\omega \in \bigoplus_{i=1}^4{\mathbb Z}\omega_i}$, 
we have $\overline{D}(S(\omega))=I(\omega,\omega)S(\omega)$. 
\end{enumerate}
\end{prop}
\begin{pf}
For (i), (ii), we refer Kac \cite{Kac}. 
For the proof of (iii), 
we use the fact that 
the initial term of 
$ch_{L(\Lambda_i)}$ is the character of the corresponding highest weight module 
of the finite dimensional Lie algebra, 
which is known to $S(0)$ ($i=0$) and $S(\omega_i)$ $(i=1,3,4)$. 
Thus we have (iii). 
For the proof of (iv), we use the fact that 
the automorphism of $D_4$ root system induces the cyclic group action 
on $\Lambda_1,\Lambda_3,\Lambda_4$. 
Then we have (iv). 
For the proof of (v), we use the theory of string function 
developped in Kac-Peterson \cite{Chevalley2}. 
Then $\chi_{\Lambda_i}$ is a product of $\eta^{-4}$ and theta function. 
Since the theta function is annihilated by $D$, 
we have (v). 
For the proof of (vi), we observe that $W_{\mathrm{cl}}$ preserves $I$. 
Then we have the result. 
\qed\end{pf}

\begin{prop}\label{631.004}
Put 
\begin{eqnarray*}
\hat{s}_0&:=&(-6)(\eta/\omega)^{-4}
\det
\bp 
\hat{\chi}_{\Lambda_1}+\hat{\chi}_{\Lambda_3}+\hat{\chi}_{\Lambda_4}&\hat{\chi}_{\Lambda_0}\\
\delta_q(\hat{\chi}_{\Lambda_1}^{\q}+\hat{\chi}_{\Lambda_3}^{\q}+\hat{\chi}_{\Lambda_4}^{\q})
&
\delta_q(\hat{\chi}_{\Lambda_0}^{\q})
\ep, \\
\hat{s}_1&:=&(\eta/\omega)^{-4}
\det
\bp 
\hat{\chi}_{\Lambda_1}+\hat{\chi}_{\Lambda_3}+\hat{\chi}_{\Lambda_4}&\hat{\chi}_{\Lambda_0}\\
\hat{\chi}_{\Lambda_1}^{\q}+\hat{\chi}_{\Lambda_3}^{\q}+\hat{\chi}_{\Lambda_4}^{\q}
&
\hat{\chi}_{\Lambda_0}^{\q}
\ep,\\
\hat{s}_2&:=&(\eta/\omega)^{-8}
(-2\hat{\chi}_{\Lambda_1}+\hat{\chi}_{\Lambda_3}+\hat{\chi}_{\Lambda_4}),\\
\hat{s}_3&:=&(\eta/\omega)^{-8}
(\hat{\chi}_{\Lambda_3}-\hat{\chi}_{\Lambda_4}).
\end{eqnarray*}
Then these are Jacobi forms 
with 
$\hat{s}_0 \in J_{0,1}$, 
$\hat{s}_1 \in J_{-2,1}$, 
$\hat{s}_2 \in J_{-4,1}$, 
$\hat{s}_3 \in J_{-4,1}$ and 
their initial terms are 
\begin{eqnarray*}
(\varphi_2)^*(\hat{s}_0)&=&
e^{\Lambda_0}\omega^{0}(S(\omega_1)+S(\omega_3)+S(\omega_4)+48), \\
(\varphi_2)^*(\hat{s}_1)&=&
e^{\Lambda_0}\omega^{4}(S(\omega_1)+S(\omega_3)+S(\omega_4)-24), \\
(\varphi_2)^*(\hat{s}_2)&=&
e^{\Lambda_0}\omega^{8}(-2S(\omega_1)+S(\omega_3)+S(\omega_4)), \\
(\varphi_2)^*(\hat{s}_3)&=&
e^{\Lambda_0}\omega^{8}(S(\omega_3)-S(\omega_4)). 
\end{eqnarray*}
\end{prop}
\begin{pf}
By Proposition \ref{631.001} (ii), 
$1$-dimensional vector spaces 
$\CC (-2\hat{\chi}_{\Lambda_1}+\hat{\chi}_{\Lambda_3}+\hat{\chi}_{\Lambda_4})$, 
$\CC (\hat{\chi}_{\Lambda_3}-\hat{\chi}_{\Lambda_4})$ 
and a $2$-dimensional vector space
$\CC (\hat{\chi}_{\Lambda_1}+\hat{\chi}_{\Lambda_3}+\hat{\chi}_{\Lambda_4})
\oplus \CC (\hat{\chi}_{\Lambda_0})$ have the $Mp_2(\ZZ)$ action. 
Then by the determinant construction (Proposition \ref{622.001}), 
we have Jacobi forms. 
Their initial terms are calculated by Proposition \ref{631.001} (iii). 
\qed\end{pf}

\begin{prop}\label{631.005}
Put 
$$
\hat{s}_4:=
\frac{1}{24}(3\mathcal{I}(\hat{s}_3,\hat{s}_3)+\mathcal{I}(\hat{s}_2,\hat{s}_2)). 
$$
Then $\hat{s}_4$ is a Jacobi form with $\hat{s}_4 \in J_{-6,2}$ 
and its initial term is 
\begin{eqnarray*}
(\varphi_2)^*(\hat{s}_4)&=&
\frac{-1}{36}
e^{2\Lambda_0}\omega^{12}
[
2S(\omega_1)^2+
2S(\omega_3)^2+
2S(\omega_4)^2+
S(\omega_1)S(\omega_3)+
S(\omega_1)S(\omega_4)+\\
&&S(\omega_3)S(\omega_4)+
24(S(\omega_1)+S(\omega_3)+S(\omega_4))
-36S(\omega_2)-288)]. 
\end{eqnarray*}
\end{prop}

\begin{pf}
By the diagram (\ref{621.003}), we have 
$\hat{s}_4 \in J_{-6,2}$. 
By the diagram (\ref{621.003}), the calculation 
of the initial terms reduces to the $\overline{D}$ action 
on $\CC[S(\omega_1),\cdots,S(\omega_4)]$. 
The action of the operator $\overline{D}$ 
on the monomials of $S(\omega_i)$ could be calculated 
by Proposition \ref{631.001} (v) and Proposition \ref{631.010}. 
\qed\end{pf}

\begin{prop}\label{631.010}
For $\omega \in \ZZ_{\geq 0}\omega_1+\cdots+\ZZ_{\geq 0}\omega_l$, 
the vector space 
$$
V:=\langle S(\omega')\,|\,\omega' \prec \omega \rangle
$$
has the following $2$ bases:
\begin{enumerate}
\item $\{S(\omega')\,|\,\omega' \prec \omega\}$, 
\item $\{\prod_{i=1}^lS(\omega_i)^{n_i}\,|\,\sum_{i=1}^ln_i\omega_i \prec \omega\}$, 
\end{enumerate}
where $\prec$ is a partial ordering on $F_{\mathrm{cl}}$ defined 
by the positive roots. 
We also have the algorithm to obtain the explicit relation 
of these bases. 
\end{prop}
\begin{pf}
The former assersion is obtained by the decomposition 
$$
S(\omega')S(\omega'')=\sum_{\omega'''\prec \omega'+\omega''} a_{\omega'''}S(\omega''')
$$
for $a_{\omega'''}\in \ZZ_{\geq 0}$ which is an easy consequence of 
Bourbaki \cite[Ch 6. \S 3]{Bourbaki}. 
The latter assersion is obtained by representation theory (cf. \cite{Koike}). 
\qed\end{pf}
\begin{prop}\label{631.006}
\begin{enumerate}
\item $\hat{s}_0,\cdots,\hat{s}_4$ are polynomial generators 
of $J_{*,*}$ over $M_*$. 
\item For $s_0,\cdots,s_4$ defined in 
(\ref{520.001}), $s_i=L^{-1}_{k_i}(\varphi_1^*(\hat{s}_i))$. In particular 
$s_i$ are polynomial generators of $S^W$ over $F(\HH) $. 
\end{enumerate}
\end{prop}
\begin{pf}
We give the proof of (i). 
By the structure theorem of Jacobi forms (Theorem \ref{623.001}), 
we should check that 
$\hat{s}_0 \in J_{0,1}\setminus \hat{E}_4J_{-4,1}$, 
$\hat{s}_1 \in J_{-2,1}\setminus \{0\}$, 
$\CC \hat{s}_2 \oplus \CC \hat{s}_3=J_{-4,1}$, 
$\hat{s}_4 \in J_{-6,2}\setminus (J_{-2,1}\cdot J_{-4,1})$. 
These are checked by the initial terms 
given by Proposition \ref{631.004} and Proposition \ref{631.005}. 
We give the proof of (ii). 
For $i=0,1,2,3$, we have $s_i=L^{-1}_{k_i}(\varphi_1^*(\hat{s}_i))$. 
For $i=4$, we have $s_i=L^{-1}_{k_i}(\varphi_1^*(\hat{s}_i))$ 
by (\ref{614.003}). 
By (i) and Theorem \ref{623.001} (ii), $L^{-1}_{k_i}(\varphi_1^*(\hat{s}_i))$ 
are polynomial generators of $S^W$ over $F(\HH) $, 
then we have (ii). 
\qed\end{pf}

We put 
$S_1:=F(\HH) [\tilde{s}_0,\cdots,\tilde{s}_3]
\oplus \tilde{s}_4 F(\HH) [\tilde{s}_0,\cdots,\tilde{s}_3]$, 
$S_2:=\tilde{s}_4^2 S^W$. 
Then we have $S^W=S_1 \oplus S_2$. 
\begin{prop}\label{631.007}
In the decomposition $S^W=S_1\oplus S_2$, 
$S_2$-component of $I^*(d\tilde{s}_i,d\tilde{s}_j) \in S^W$ 
is $0$ for $i,j=-1,\cdots,4$. 
\end{prop}

\begin{pf}
For the cases $i=-1$ or $j=-1$, the assertions are shown by (\ref{240.002}). 
We assume $i,j \geq 0$. 
We put $J_1:=M_*[\hat{s}_0,\cdots,\hat{s}_3]
\oplus \hat{s}_4M_*[\hat{s}_0,\cdots,\hat{s}_3]$, 
$J_2:=\hat{s}_4^2 J_{*,*}$. Then we have $J_{*,*}=J_1\oplus J_2$. 
By (\ref{614.003}), we should only prove that 
$J_2$-component of the Jacobi form $\I(\hat{s}_i,\hat{s}_j)$ is $0$. 
We show it by checking its weight and index. 
The weight and index of this Jacobi form 
is $k_i+k_j+2$ and $m_i+m_j$. 
However the $J_2$-component of such weight and index must be $0$ 
because of the structure theorem of Wirthm\"uller. 
Thus we have the result. 
\qed\end{pf}

We put $J_{k,m}^{\mathrm{cusp}}:=\mathrm{ker}(J_{k,m} \to J^0_{k,m})$. 
Then we have $J_{k,m}^{\mathrm{cusp}}=(\hat{E}_4^3-\hat{E}_6^2)J_{k-12,m}$. 
Then we have the following Lemma. 
\begin{lem}\label{623.002}
If $J_{k-12,m}=\{0\}$, then 
$J_{k,m} \to J^0_{k,m}$ is an isomorphism. 
\end{lem}

\begin{prop}\label{631.008}
We have  
\begin{eqnarray*}
\I(\hat{s}_0,\hat{s}_3)&=&-\frac{1}{6}(3\hat{E}_4 \hat{s}_1+\hat{E}_6\hat{s}_2)\hat{s}_3,\\
\I(\hat{s}_1,\hat{s}_3)&=&-\frac{1}{6} (2\hat{s}_0+\hat{E}_4 \hat{s}_2)\hat{s}_3,\\
\I(\hat{s}_2,\hat{s}_3)&=&-\frac{1}{3} \hat{s}_1 \hat{s}_3,\\
\I(\hat{s}_3,\hat{s}_3)&=&4\hat{s}_4-\frac{1}{9} \hat{s}_1 \hat{s}_2\\
\I(\hat{s}_2,\hat{s}_2)&=&12\hat{s}_4+\frac{1}{3} \hat{s}_1\hat{s}_2,\\
\I(\hat{s}_1,\hat{s}_2)&=&-\frac{1}{3} \hat{s}_0\hat{s}_2-\frac{1}{4}
 \hat{E}_4 \hat{s}_3 \hat{s}_3+\frac{1}{12} \hat{E}_4 \hat{s}_2 \hat{s}_2,\\
\I(\hat{s}_0,\hat{s}_2)&=&-\frac{1}{4} \hat{E}_6 \hat{s}_3 \hat{s}_3-
\frac{1}{12}(6\hat{E}_4 \hat{s}_1-\hat{E}_6 \hat{s}_2)\hat{s}_2,\\
\I(\hat{s}_1,\hat{s}_1)&=&6 \hat{E}_4 \hat{s}_4-\frac{1}{6} \hat{s}_0\hat{s}_1
-\frac{1}{24} \hat{E}_6(\hat{s}_2\hat{s}_2+3\hat{s}_3\hat{s}_3),\\
\I(\hat{s}_0,\hat{s}_1)&=&6 \hat{E}_6 \hat{s}_4
-\frac{1}{3} \hat{E}_4 \hat{s}_1\hat{s}_1-\frac{1}{24} \hat{E}_4^2(\hat{s}_2\hat{s}_2+3\hat{s}_3\hat{s}_3),\\
\I(\hat{s}_0,\hat{s}_0)&=&6 \hat{E}_4^2 \hat{s}_4
-\frac{1}{3} \hat{E}_6 \hat{s}_1\hat{s}_1-\frac{1}{6} \hat{E}_4 \hat{s}_0\hat{s}_1
-\frac{1}{24} 
\hat{E}_4\hat{E}_6(\hat{s}_2\hat{s}_2+3\hat{s}_3\hat{s}_3),\\
\I(\hat{s}_3,\hat{s}_4)&=&\frac{1}{432} \hat{s}_3(8\hat{s}_1 \hat{s}_1
+8\hat{s}_0\hat{s}_2+\hat{E}_4 \hat{s}_2\hat{s}_2+3 \hat{E}_4 \hat{s}_3\hat{s}_3),\\
\I(\hat{s}_2,\hat{s}_4)&=&\frac{1}{36}\hat{s}_0\hat{s}_3\hat{s}_3
+\frac{1}{54} \hat{s}_1\hat{s}_1\hat{s}_2
-\frac{1}{108} \hat{s}_0\hat{s}_2\hat{s}_2
+\frac{1}{432}\hat{E}_4\hat{s}_2(\hat{s}_2\hat{s}_2+3\hat{s}_3\hat{s}_3),\\
\I(\hat{s}_1,\hat{s}_4)&=&-\frac{1}{2}\hat{s}_0\hat{s}_4
+\frac{1}{144} \hat{E}_4\hat{s}_1(\hat{s}_2\hat{s}_2+3\hat{s}_3\hat{s}_3)
+\frac{1}{864} \hat{E}_6\hat{s}_2(\hat{s}_2+3\hat{s}_3)(-\hat{s}_2+3\hat{s}_3),\\
\I(\hat{s}_0,\hat{s}_4)&=&-\frac{5}{6} \hat{E}_4\hat{s}_1\hat{s}_4
+\frac{1}{144} \hat{E}_6\hat{s}_1(\hat{s}_2\hat{s}_2+3\hat{s}_3\hat{s}_3)
+\frac{1}{864} \hat{E}_4^2 \hat{s}_2(\hat{s}_2+3\hat{s}_3)(-\hat{s}_2+3\hat{s}_3),\\
\I(\hat{s}_4,\hat{s}_4)&=&\frac{1}{432}\hat{s}_4(8\hat{s}_1^2
+\hat{E}_4 \hat{s}_2^2+3 \hat{E}_4 \hat{s}_3^2)
-\frac{5}{7776} \hat{s}_0\hat{s}_1(\hat{s}_2^2+3\hat{s}_3^2)
-\frac{1}{31104} \hat{E}_6 \hat{s}_2^4
\\
&&-\frac{1}{5184} \hat{E}_4 \hat{s}_1\hat{s}_2(\hat{s}_2+3\hat{s}_3)(-\hat{s}_2+3\hat{s}_3)
-\frac{1}{10368} \hat{E}_6 \hat{s}_3^2(2\hat{s}_2^2+3\hat{s}_3^2),\\
\mathcal{D}(\hat{s}_1)&=&-\frac{1}{3}\hat{s}_0,\\
\mathcal{D}(\hat{s}_0)&=&-\frac{2}{3}\hat{E}_4 \hat{s}_1.
\end{eqnarray*}
\end{prop}
\begin{pf}
For each equation, the both sides are Jacobi forms 
of same weight and index. 
For these weight and index, 
we could check that the space of cusp forms 
are $\{0\}$ by Lemma \ref{623.002}. 
Thus we should only check 
that the initial terms of both sides coincide. 
This could be done by 
calculating the $\overline{D}$ action 
on $\CC[S(\omega_1),\cdots,S(\omega_4)]$. 
which is discussed in Proposition \ref{631.010}. 
\qed\end{pf}
\section{Differential relations satisfied by modular forms}
In Section \ref{170831.2}, 
we show that $\eta^4\chi_{\Lambda_i}^{\q}$ $(i=0,1)$ 
satisfy the linear differential equation called ``the Kaneko-Zagier equation" 
in Proposition \ref{641.002} eq. (\ref{170831.4}). 
It is used to give the flat generator system 
in Proposition \ref{170831.33} (ii).
In Section \ref{170831.3}, 
we show that $\eta^4\chi_{\Lambda_i}^{\q}$ $(i=0,1)$ 
satisfy the non-linear differential equations called ``the Halphen's equations"
in Proposition \ref{641.003} eq. (\ref{170831.1}). 
It is used to give a description of the potential of the Frobenius structure 
in Proposition \ref{170831.33}. 

\subsection{The Kaneko-Zagier equations}\label{170831.2}
For $k \in \frac{1}{2}\mathbb{Z}$, we consider the following differential equation:
\begin{equation}
\p_{k+2}\p_k f=\frac{k}{12}\frac{k+2}{12}E_4f, 
\end{equation}
which we call ``the Kaneko-Zagier equation" \cite{Kaneko-Koike} 
for $k \in \frac{1}{2}\mathbb{Z}$, 
where we denote $f'-\frac{k}{12}E_2f$ by $\p_kf$ for $f \in F(\HH)$ 
defined in Proposition \ref{614.005}. 
\begin{prop}\label{641.002}
\begin{enumerate}
\item For $i=0,1$, $\eta^4\chi_{\Lambda_i}^{\q}$ 
satisfies the Kaneko-Zagier equation for $k=2$:
\begin{equation}\label{170831.4}
\p_4\p_2(\eta^4\chi_{\Lambda_i}^{\q})=
\frac{2}{12}\frac{2+2}{12}E_4
\eta^4 \chi_{\Lambda_i}^{\q}. 
\end{equation}
\item $\det K=4$ 
for $K:=\bp \chi_{\Lambda_0}^{\q}&\chi_{\Lambda_1}^{\q}\\
\eta^{-4}(\chi_{\Lambda_0}^{\q})'&
\eta^{-4}(\chi_{\Lambda_1}^{\q})'
\ep$.
\item The matrix $K$ satisfies the following relation:
\begin{equation}
\tr{(CK)}MCK=\bp 48&0\\0&16\ep
\end{equation}
for $C:=\bp1/6&0\\0&1\ep$, 
$M:=
\bp \tilde{E}_4^2&\tilde{E}_6\\
\tilde{E}_6&\tilde{E}_4\ep$. 
\end{enumerate}
\end{prop}
\begin{pf}
By Proposition \ref{631.008}, 
we have $\mathcal{D}(\hat{s}_1)=\frac{-1}{3}\hat{s}_0$, 
$\mathcal{D}(\hat{s}_0)=\frac{-2}{3}\hat{E}_4 \hat{s}_1$. 
Then we have 
\begin{equation}
\mathcal{D}\mathcal{D}(\hat{s}_1)=
\frac{2}{9}\hat{E}_4 \hat{s}_1. 
\end{equation}
Substituting 
$
\hat{s}_1=(\eta/\omega)^{-4}
\det
\bp 
\hat{\chi}_{\Lambda_1}+\hat{\chi}_{\Lambda_3}+\hat{\chi}_{\Lambda_4}&\hat{\chi}_{\Lambda_0}\\
3\hat{\chi}_{\Lambda_1}^{\q}
&
\hat{\chi}_{\Lambda_0}^{\q}
\ep
$
in Proposition \ref{631.004} 
(here we used $\hat{\chi}_{\Lambda_1}^{\q}=\hat{\chi}_{\Lambda_3}^{\q}=\hat{\chi}_{\Lambda_4}^{\q}$ 
which is shown in Proposition \ref{631.001}), 
and 
using 
$
\mathcal{D}(\hat{f}\hat{F})
=2\delta_q(\hat{f})E(\hat{F})
+\hat{f}\mathcal{D}(\hat{F})
$ 
in (\ref{613.001}), 
$\delta_q(\eta/\omega)=0$ 
in (\ref{613.003}) 
and 
$\mathcal{D}(\hat{\chi}_i)=0$ $(i=0,1)$ in Proposition \ref{631.001}, 
we have 
$$
\det
\bp 
\hat{\chi}_{\Lambda_1}+\hat{\chi}_{\Lambda_3}+\hat{\chi}_{\Lambda_4}&\hat{\chi}_{\Lambda_0}\\
4\delta_q\delta_q[3\hat{\chi}_{\Lambda_1}^{\q}]
&
4\delta_q\delta_q[\hat{\chi}_{\Lambda_0}^{\q}]
\ep
=
\frac{2}{9}\hat{E}_4 
\det
\bp 
\hat{\chi}_{\Lambda_1}+\hat{\chi}_{\Lambda_3}+\hat{\chi}_{\Lambda_4}&\hat{\chi}_{\Lambda_0}\\
3\hat{\chi}_{\Lambda_1}^{\q}
&
\hat{\chi}_{\Lambda_0}^{\q}
\ep. 
$$
Comparing of coefficients of 
$\hat{\chi}_{\Lambda_0}$ and $\hat{\chi}_{\Lambda_1}+\hat{\chi}_{\Lambda_3}+\hat{\chi}_{\Lambda_4}$, 
we have 
\begin{equation}
\delta_q\delta_q(\hat{\chi}_{\Lambda_i}^{\q})
=\frac{1}{18}\hat{E}_4\hat{\chi}_{\Lambda_i}^{\q}
\end{equation}
for $i=0,1$. 
Then we have 
\begin{equation}
X_{\omega}((\eta/\omega)^4\hat{\chi}_{\Lambda_i}^{\q})
=(-4)(\eta/\omega)^4\hat{\chi}_{\Lambda_i}^{\q},\ 
\delta_q\delta_q((\eta/\omega)^4\hat{\chi}_{\Lambda_i}^{\q})
=
\frac{1}{18}\hat{E}_4
(\eta/\omega)^4\hat{\chi}_{\Lambda_i}^{\q}. 
\end{equation}
This means that $(\eta/\omega)^4\chi_{\Lambda_i}^{\q}$ satisfies differential 
equation of Proposition \ref{720.002} for $k=2$. 
Then by Proposition \ref{720.002}, we have (i). 
The $q$-expansions of $\eta^4 \chi_{\Lambda_i}^{\q}$ are 
\begin{equation}
\eta^4 \chi_{\Lambda_0}^{\q}
=1+\cdots \in \CC\{q\}, 
\eta^4 \chi_{\Lambda_1}^{\q}
=8q^{\frac{1}{2}}+\cdots \in q^{\frac{1}{2}}\CC\{q\}. 
\end{equation}
Applying (\ref{710.004}), we have (ii). 
By Proposition \ref{710.003} (iii), we have (iii). 
\qed\end{pf}

\subsection{The Halphen's equations}\label{170831.3}
For $E_2,E_4,E_6$, we see that by the famous relations 
\begin{equation}\label{641.001}
(E_2)'=\frac{1}{12}(E_2^2-E_4)
,\quad
(E_4)'=\frac{1}{3}(E_2E_4-E_6)
,\quad
(E_6)'=\frac{1}{2}(E_2E_6-E_4^2).
\end{equation}
(cf. Kaneko-Koike \cite{Kaneko-Koike}), 
the $\CC$-algebra generated by $E_2,E_4,E_6$ 
has a structure of differential algebra. 

We assert that the $\CC$-algebra 
generated by 
$\eta^4 \chi_{\Lambda_0}^{\q},
\eta^4 \chi_{\Lambda_1}^{\q}
,\eta'/\eta$ 
also has a structure of differential algebra 
which contains the above differential algebra. 

For that purpose, we first introduce the following functions 
\begin{subequations}\label{641.004}
\begin{align}
\xi_2&:=2(\log \sum_{n \in \ZZ}q^{(n-1/2)^2/2})',\\
\xi_3&:=2(\log \sum_{n \in \ZZ}q^{n^2/2})',\\
\xi_4&:=2(\log \sum_{n \in \ZZ}(-)^n q^{n^2/2})'.
\end{align}
\end{subequations}
They satisfy the following differential relations 
studied by Halphen \cite{Halphen}:
\begin{subequations}\label{170831.1}
\begin{eqnarray}
\xi_2'&=&\xi_2\xi_3+\xi_2\xi_4-\xi_3\xi_4,\label{641.006}\\
\xi_3'&=&\xi_2\xi_3-\xi_2\xi_4+\xi_3\xi_4,\label{641.007}\\
\xi_4'&=&-\xi_2\xi_3+\xi_2\xi_4+\xi_3\xi_4,\label{641.008}
\end{eqnarray}
\end{subequations}
which we call ``the Halphen's equations". 
For these equations, we refer Ohyama \cite{Ohyama}. 
\begin{prop}\label{641.003}
$\eta^4 \chi_{\Lambda_0}^{\q},
\eta^4 \chi_{\Lambda_1}^{\q}
,\eta'/\eta$ satisfy the following relations. 
\begin{subequations}
\begin{align}
(\eta^4 \chi_{\Lambda_0}^{\q})'&=
4(\eta^4 \chi_{\Lambda_0}^{\q})
\frac{\eta'}{\eta}
-\frac{1}{6}(\eta^4 \chi_{\Lambda_0}^{\q})^2
+\frac{1}{2}(\eta^4 \chi_{\Lambda_1}^{\q})^2,\\
(\eta^4 \chi_{\Lambda_1}^{\q})'
&=
4(\eta^4 \chi_{\Lambda_1}^{\q})
\frac{\eta'}{\eta}
+\frac{1}{3}(\eta^4 \chi_{\Lambda_0}^{\q})(\eta^4 \chi_{\Lambda_1}^{\q}),\\
(\frac{\eta'}{\eta})'&=
2\left(\frac{\eta'}{\eta}\right)^2-\frac{1}{2^5 3^2}
[(\eta^4 \chi_{\Lambda_0}^{\q})^2
+3(\eta^4 \chi_{\Lambda_1}^{\q})^2],
\end{align}
\end{subequations}
and 
\begin{eqnarray}
E_2&=&24 \frac{\eta'}{\eta},\\
E_4&=&(\eta^4 \chi_{\Lambda_0}^{\q})^2
+3(\eta^4 \chi_{\Lambda_1}^{\q})^2,\\
E_6&=&(\eta^4 \chi_{\Lambda_0}^{\q})^3
-9(\eta^4 \chi_{\Lambda_0}^{\q})(\eta^4 \chi_{\Lambda_1}^{\q})^2.
\end{eqnarray}
\end{prop}
We show the following proposition. 
\begin{prop}\label{641.005}
\begin{eqnarray*}
\eta^4 \chi_{\Lambda_0}^{\q}&=&(4\xi_2-2\xi_3-2\xi_4),\\
\eta^4 \chi_{\Lambda_1}^{\q}&=&(2\xi_3-2\xi_4),\\
\frac{\eta'}{\eta}&=&\frac{1}{6}(\xi_2+\xi_3+\xi_4).
\end{eqnarray*}
\end{prop}
\begin{pf}
We first remark that 
$\frac{\eta'}{\eta}=\frac{1}{6}(\xi_2+\xi_3+\xi_4)$ 
by Jacobi's derivative formula (cf. \cite{Mumford}). 
For $i=2,3,4$, we have 
$$
\p_{2+2}\p_{2}\xi_i-\frac{2}{12}\frac{2+2}{12}E_4\xi_i
=2\xi_2\xi_3\xi_4
$$
for $\p_k f=f'-\frac{k}{12}E_2 f$ 
since 
\begin{eqnarray*}
&&\p_{2+2}\p_{2}\xi_i-\frac{2}{12}\frac{2+2}{12}E_4\xi_i\\
&=&
(\xi_i)''-\frac{2+1}{6}E_2 (\xi_i)'+
\frac{2(2+1)}{12}E_2'\xi_i\\
&=&
(\xi_i)''-\frac{2+1}{6}(4(\xi_2+\xi_3+\xi_4))(\xi_i)'+
\frac{2(2+1)}{12}(4(\xi_2+\xi_3+\xi_4))'\xi_i\\
&=&
2\xi_2\xi_3\xi_4.
\end{eqnarray*}
Here we use 
$
E_2'=\frac{1}{12}(E_2^2-E_4)
$
for 1st equality, 
use $E_2=4(\xi_2+\xi_3+\xi_4)$ for 2nd equality 
and use (\ref{170831.1}) for the last equality. 

Thus if $a+b+c=0$ for $a,b,c \in \CC$, then $a\xi_2+b\xi_3+c\xi_4$ 
satisfies the Kaneko-Zagier equation for $k=2$. 

Since $\eta^4 \chi_{\Lambda_0}^{\q},\eta^4 \chi_{\Lambda_1}^{\q}$ also satisfy the 
Kaneko-Zagier equation for $k=2$, 
we obtain the result by comparing the leading terms 
of the $q$-expansions. 
\qed\end{pf}

\begin{pf}[of Proposition \ref{641.003}]
In Ohyama \cite{Ohyama}, the following relations are obtained:
\begin{equation}\label{641.009}
\tilde{h}_2=-1/48 E_4,
\quad
\tilde{h}_3=1/2^53^3 E_6,
\end{equation}
where $h_1:=\xi_2+\xi_3+\xi_4$, $h_2:=\xi_2\xi_3+\xi_2\xi_4+\xi_3\xi_4$, $h_3:=\xi_2\xi_3\xi_4$, 
$\tilde{h}_2:=h_2-1/3h_1^2$, 
$\tilde{h}_3:=h_3-1/3h_1h_2+2/27h_1^3$. 
By Proposition \ref{641.005}, the equations 
(\ref{641.006})--(\ref{641.008}) and (\ref{641.009}), 
we have the results. 
\qed\end{pf}

\section{Duality for the Kaneko-Zagier equation}
Let $V_k$ be a solution space of the Kaneko-Zagier equation 
for $k \in \frac{1}{2}\mathbb{Z}$. 
For $k \neq 0,4$, we find the duality between $V_{k}$ and $V_{4-k}$, 
which appeared already in Proposition \ref{170831.33} (i) eq. (\ref{520.012}) 
for $k=2$ and \cite[(3.25)]{E6} for $k=3$. 
This duality is important for the study of the differential algebra 
generated by $V_k$ which will be used to the explicit construction 
of the potential of the Frobenius structure for the elliptic root 
systems of type $E_6^{(1,1)}$ ($k=3$ case) and $E_7^{(1,1)}$ ($k=7/2$ case). 

In this section, we first formulate the Kaneko-Zagier equations 
as a connection on $\CC^* \times \HH$ (eq. (\ref{170831.11})). 
Then we have a duality for any $k \in \frac{1}{2}\mathbb{Z}$ (eq. (\ref{170831.12})). 
This formulation is analogous to the formulation of Saito \cite[\S 5.4]{period}. 
Then we restrict ourselves to the cases of $0 < k < 4$ 
and give the duality for the Kaneko-Zagier equations (eq. (\ref{170831.13})). 

We put $M:=\CC^* \times \HH$. 
We denote by $\Oo_M$ the sheaf of holomorphic functions on $M$ 
and by $\Omega_M$ the sheaf of holomorphic $1$-forms on $M$. 
Let 
$$
X_1:=\omega\frac{\p}{\p \omega}, \ X_2:=\delta_q
$$
be vector fields which give a frame on each $p \in M$. 
Let 
$$
\zeta_1:=
\frac{d\omega}{\omega}-\frac{1}{24}\hat{E}_2\frac{dq}{q},\quad
\zeta_2:=\omega^4\frac{dq}{q}
$$
be dual 1-forms on $M$. 

We define the non-degenerate $\Oo_M$-symmetric bilinear form 
$h:\Omega_M \times \Omega_M \to \Oo_M$ by 
$$
\bp 
h(\zeta_1,\zeta_1)&
h(\zeta_1,\zeta_2)\\
h(\zeta_2,\zeta_1)&
h(\zeta_2,\zeta_2)
\ep
=\frac{1}{(\frac{\eta}{\omega})^{24}}
\bp \frac{-1}{24}&0\\0&1\ep
\bp \hat{E}_4^2&\hat{E}_6\\
\hat{E}_6&\hat{E}_4\ep
\bp \frac{-1}{24}&0\\0&1\ep. 
$$

For $k \in \frac{1}{2}\ZZ$, 
we define a connection $\nabla^{(k)}:\Omega_M \to 
\Omega_M\otimes_{\Oo_M} \Omega_M$ by 
\begin{equation}\label{170831.11}
(\nabla^{(k)} \zeta_1,\nabla^{(k)} \zeta_2)
=(\zeta_1,\zeta_2)\Gamma^{(k)}, 
\end{equation}
where the connection form 
$\Gamma^{(k)}:=\Gamma^{(k)}_1\zeta_1+\Gamma^{(k)}_2\zeta_2$ 
is given by 
$$
\Gamma^{(k)}_1
=
\bp 2k&0\\
0&2(k+2)\ep
,\quad
\Gamma^{(k)}_2
=
\bp 0&2k\\
\frac{1}{24}\frac{k+2}{12}\hat{E}_4&0\ep. 
$$
We see that it is a flat and torsion free connection on $\Omega_M$. 

This connection has the following property which 
is obtained by direct calculation.
\begin{prop}\label{720.001}
We have 
\begin{equation}\label{170831.12}
dh(\omega,\omega')=
h(\nabla^{(k)}\omega,\omega')+h(\omega,\nabla^{(4-k)}\omega')
\end{equation}
for local sections $\omega,\omega'$ of $\Omega_M$. 
\end{prop}

By torsion freeness of $\nabla^{(k)}$, the flat section is a closed form. 
Thus the differential equations 
$\nabla^{(k)}dF=0$ have 3-dimensional solutions which 
contain constant functions. 

Hereafter we also assume the homogeneity. Then we have the following differential equations. \\

\begin{prop}\label{720.002}
For $k \in \frac{1}{2}\ZZ$ and for a local section $F$ of $\Oo_M$, 
the following conditions are equivalent:
\begin{eqnarray*}
&&(1)\ X_1 F=(-2k)F,\quad
\nabla^{(k)}dF=0.\\
&&(2)\ X_1 F=(-2k)F,\quad
\delta_q \delta_q F=\frac{k}{12}\frac{k+2}{12}\hat{E}_4F.\\
&&(3) \hbox{There exists a local section $f$ of $\Oo_{\HH}$ s.t. 
$F=\omega^{-2k}f$ and $\p_{k+2}\p_k f=\frac{k}{12}\frac{k+2}{12}E_4f$}.
\end{eqnarray*}
\end{prop}
For each $k \in \frac{1}{2}\ZZ$, we have 2-dimensional solutions. 
They contain non-zero constant functions iff $k=0$. 
Thus these solutions have pairing induced by $h$ which is non-degenerate 
for $k\neq 0,4$, and degenerate for $k=0,4$. 

Hereafter we assume $0<k<4$. 
For these cases, we shall write down the pairing explicitly. 
We have one solution $f^{(k)}_1 \in \CC\{\{q\}\}$ 
whose leading term is $1$. 
We also have another solution $f^{(k)}_2 \in q^{\alpha}\CC\{\{q\}\}$ 
whose leading term is $q^{\alpha}$ where $\alpha=\frac{k+1}{6}$. 
By these conditions, $f^{(k)}_j$ ($j=1,2$) are uniquely determined. 
Their $q$-expansions are
\begin{eqnarray*}
f^{(k)}_1&=&1+a_1q+\cdots,\\
f^{(k)}_2&=&q^{\alpha}+b_1q^{\alpha+1}+\cdots, 
\end{eqnarray*}
with $a_1=\frac{12k(k+1)}{5-k}$, 
$b_1=\frac{4(k+1)(2k-1)}{k+7}$. 

We see that 
\begin{equation}\label{710.004}
\det
\bp 
f^{(k)}_1&
f^{(k)}_2\\
\p_k f^{(k)}_1&
\p_k f^{(k)}_2
\ep
=\alpha\eta^{24\alpha},
\end{equation}
by Wronskian and $q$-expansion. 

We put 
$F^{(k)}_1=\omega^{-2k}f^{(k)}_1$, 
$F^{(k)}_2=\omega^{-2k}f^{(k)}_2$. 

\begin{prop}\label{710.003}
For $0 < k < 4$, we have 
\begin{enumerate}
\item 
$$
\bp 
h(dF^{(k)}_1,dF^{(4-k)}_1)&
h(dF^{(k)}_1,dF^{(4-k)}_2)\\
h(dF^{(k)}_2,dF^{(4-k)}_1)&
h(dF^{(k)}_2,dF^{(4-k)}_2)
\ep
=
\bp 12k(4-k)&0\\
0&\frac{(k+1)(5-k)}{36}\ep.
$$
\item 
Put 
$$
\mathcal{F}^{(k)}:=
\bp 
(-2k)f^{(k)}_1&(-2k)f^{(k)}_2\\
\p_{k}f^{(k)}_1
&
\p_{k}f^{(k)}_2
\ep. 
$$
Then we have a duality:
\begin{equation}\label{170831.13}
\tr{\mathcal{F}^{(k)}}
\frac{1}{\eta^{24}}
\bp \frac{-1}{24}&0\\0&1\ep
\bp E_4^2&E_6\\
E_6&E_4\ep
\bp \frac{-1}{24}&0\\0&1\ep
\mathcal{F}^{(4-k)}
=
\bp 12k(4-k)&0\\
0&\frac{(k+1)(5-k)}{36}\ep. 
\end{equation}
\end{enumerate}
\end{prop}
\begin{pf}
Put 
$$
\hat{\mathcal{F}}^{(k)}
:=
\bp 
(-2k)\omega^{-2k}f^{(k)}_1&(-2k)\omega^{-2k}f^{(k)}_2\\
\omega^{-4-2k}\p_{k}f^{(k)}_1
&
\omega^{-4-2k}\p_{k}f^{(k)}_2
\ep. 
$$
Then 
\begin{eqnarray*}
&&
\bp 
h(dF^{(k)}_1,dF^{(4-k)}_1)&
h(dF^{(k)}_1,dF^{(4-k)}_2)\\
h(dF^{(k)}_2,dF^{(4-k)}_1)&
h(dF^{(k)}_2,dF^{(4-k)}_2)
\ep\\
=&&
\bp 
X_1F^{(k)}_1&
X_2F^{(k)}_1\\
X_1F^{(k)}_2&
X_2F^{(k)}_2
\ep
\bp 
h(\zeta_1,\zeta_1)&
h(\zeta_1,\zeta_2)\\
h(\zeta_2,\zeta_1)&
h(\zeta_2,\zeta_2)
\ep
\bp 
X_1F^{(4-k)}_1&X_1F^{(4-k)}_2\\
X_2F^{(4-k)}_1
&
X_2F^{(4-k)}_2
\ep\\
=&&
\tr{\hat{\mathcal{F}}^{(k)}}
\bp 
h(\zeta_1,\zeta_1)&
h(\zeta_1,\zeta_2)\\
h(\zeta_2,\zeta_1)&
h(\zeta_2,\zeta_2)
\ep
\hat{\mathcal{F}}^{(4-k)}
\end{eqnarray*}
must be constant matrix 
and the $q$-expansion of the last expression could be calculated 
by 
$$
\bp 
h(\zeta_1,\zeta_1)&
h(\zeta_1,\zeta_2)\\
h(\zeta_2,\zeta_1)&
h(\zeta_2,\zeta_2)
\ep
$$
and the $q$-expansions of $f^{(k)}_j$ $(j=1,2)$, we obtain (i). 

Since 
$$
\lim_{\omega \to 1}
\hat{\mathcal{F}}^{(k)}
=\mathcal{F}^{(k)}, 
$$
we obtain (ii) by the proof of (i). 
\qed\end{pf}

\appendix
\section{A Lemma on the connection}
In this appendix, we give a lemma on the connection which was used in 
the proof of Proposition \ref{170831.32}. 

The idea is as follows. 
If we have a vector bundle $E$ on $X$ and 
the flat connection $\nabla$ on $E$ 
and we have a submersion $\varphi:X \to Y$ with 
each fiber is contractible. 
Then the flat frame of $E$ on $X$ could be obtained 
by 2 steps: 
(1) construct frame $v_1,\cdots,v_l$ of $E$ which are  
flat along every fiber of $\varphi$ 
(2) find functions $f_1,\cdots,f_l$ on $Y$ 
s.t. $\sum_{i=1}^l f_iv_i$ is flat on $X$. 

The following lemma enables us to realize step (1) 
in the setting of Proposition \ref{170831.32}. 
We denote $S^W$ by $S$ and $F(\HH)$ by $R$ 
in this appendix. 
\begin{lem}\label{710.001}
Let $E_1 \subset E_2 \subset E_3$ 
are graded $S$-free modules. 
The module $E_3$ has a flat connection 
$$
\nabla:E_3 \to \Omega_{S/\CC}\otimes_{S}E_3
$$
of degree $0$, preserving submodules $E_2,E_1$ 
respectively. 
We assume that 
$u_1,\cdots,u_l$ are $S$-free basis of $E_1$, 
each $u_i$ is homogeneous of degree less than $r$. 
$u_1,\cdots,u_l,v_1,\cdots,v_m$ are 
$S$-free basis of $E_2$, 
each $v_i$ is homogeneous of degree $r$. 
$u_1,\cdots,u_l,v_1,\cdots,v_m,w_1,\cdots,w_n$ 
 are $S$-free basis of $E_3$, 
each $w_i$ is homogeneous of degree greater than $r$, 
for some $r \in \ZZ_{\geq 0}$. 
We also assume that 
\begin{eqnarray*}
\nabla u_i&=&
\sum_{j=1}^{l}a_{ij}\otimes u_j\ (1 \leq i \leq l),\\
\nabla v_i&=&
\sum_{j=1}^{m}b_{ij}\otimes v_j+\sum_{j=1}^{l}c_{ij}u_j,\ (1 \leq i \leq m),\\
\nabla w_i&=&
\sum_{j=1}^{n}e^1_{ij}\otimes w_j+\sum_{j=1}^{m}e^2_{ij}v_j
+\sum_{j=1}^{l}e^3_{ij}\otimes u_j\ (1 \leq i \leq n),
\end{eqnarray*}
for $a_{ij} \in \Omega_{R/\CC},b_{ij},c_{ij},e^k_{ij}\in \Omega_{S/\CC}$. 
\begin{enumerate}
\item For $1 \leq i,j \leq m$, $b_{ij} \in \Omega_{R/\CC}$. 
\item For $1 \leq i \leq m,\ 1 \leq j \leq l$, 
there exists uniquely $d_{ij} \in S$ which is homogeneous and 
$d_{S/R}(d_{ij})=\pi(c_{ij})$ for $\pi:\Omega_{S/\CC} \to \Omega_{S/R}$. 
\item For $1 \leq i \leq m$, 
put $\tilde{v}_i:=v_i-\sum_{j=1}^ld_{ij}u_j$. 
Then $\tilde{v}_i$ is homogeneous of degree $k$ and 
$u_i\,(1 \leq i \leq l),\,\tilde{v}_i\,(1 \leq i \leq m),\,
w_i\,(1 \leq i \leq n)$ are $S$-free basis of $E_3$. 
By this basis, $\nabla$ is represented as 
\begin{eqnarray*}
\nabla u_i&=&
\sum_{j=1}^{l}a_{ij}\otimes u_j\ (1 \leq i \leq l),\\
\nabla \tilde{v}_i&=&
\sum_{j=1}^{m}b_{ij}\otimes \tilde{v}_j+\sum_{j=1}^{l}c_{ij}u_j,\ (1 \leq i \leq m),\\
\nabla w_i&=&
\sum_{j=1}^{n}e^1_{ij}\otimes w_j+\sum_{j=1}^{m}e^2_{ij}
[\tilde{v}_j+\sum_{k=1}^{l}d_{jk}u_k]
+\sum_{j=1}^{l}e^3_{ij}\otimes u_j\ (1 \leq i \leq n). 
\end{eqnarray*}
\end{enumerate}
\end{lem}
\begin{pf}
For the proof of (i), we obtain it by the degree condition. 
We give the proof of (ii). 
Let 
$$
\overline{\nabla}:E_2 \to \Omega_{S/R}\otimes_{S}E_2
$$
be the composed mapping of the connection on $E_2$:
$\nabla:E_2 \to \Omega_{S/\CC}\otimes_S E_2$ 
and the mapping $\Omega_{S/\CC}\otimes_S E_2 
\to \Omega_{S/R}\otimes_S E_2$ induced by 
$\pi:\Omega_{S/\CC}\to \Omega_{S/R}$. 
Then $\overline{\nabla}$ is also connection and 
it is flat. 
By the assumption $a_{ij} \in \Omega_{R/\CC}$ $(1 \leq i,j \leq l)$, 
we have $\overline{\nabla}u_i=0$ $(1 \leq i \leq l)$. 
By (i), we have $\overline{\nabla}v_i=\sum_{j=1}^{l}\pi(c_{ij})\otimes u_j$ 
($1 \leq i \leq m$). 
By the flatness of $\overline{\nabla}$ and 
$\overline{\nabla}u_i=0$, we have $d_{S/R}(\pi(c_{ij}))=0$. 

Since $\pi(c_{ij})$ is homogeneous of degree $>0$ 
and the following sequence 
$$
\begin{CD}
0 @>>> R @>>> S @>{d_{S/R}}>> \Omega_{S/R} @>{d_{S/R}}>> \Omega^2_{S/R}
\end{CD}
$$
is exact which will be shown in Lemma \ref{710.002}, 
we have the unique $d_{ij} \in S$ s.t. 
$d_{S/R}(d_{ij})=\pi(c_{ij})$. 

We give the proof of (iii). 
By the construction of $\tilde{v}_i$ ($1 \leq i \leq m$) and 
$\overline{\nabla}\tilde{v}_i=0$, we have 
$$
\nabla \tilde{v}_i
=\sum_{j=1}^{m}b_{ij}\otimes \tilde{v}_j+\sum_{j=1}^lf_{ij}\otimes u_j
$$
for some $f_{ij} \in S \otimes_R \Omega_{R/\CC}$. 
Since $\nabla$ is flat and $\Omega^2_{R/\CC}=0$, 
the entries of 
the connection matrix is closed. Thus we have $d_{S/\CC}f_{ij}=0$. 
Then $f_{ij}$ must be an element of $\Omega_{R/\CC}$. 
By the degree condition, $f_{ij}$ must be $0$. 
The other part is a direct consequence of this result. 
\qed\end{pf}

\begin{lem}\label{710.002} 
The following sequence 
$$
\begin{CD}
0 @>>> R @>>> S @>{d_{S/R}}>> \Omega_{S/R} @>{d_{S/R}}>> \Omega^2_{S/R}
\end{CD}
$$
is exact and 
if the homogeneous element 
$\omega \in \Omega_{S/R}$ 
satisfies $d_{S/R}\omega=0$, then 
we have an algorithm to 
obtain $\eta \in S$ 
s.t. $d_{S/R}(\eta)=\omega$. 
\end{lem}
\begin{pf}
Since $S$ is a polynomial algebra over $R$, 
the exactness of the sequence reduces to the one for the case of 
a polynomial algebra over $\QQ$. 
For this case, the exactness of the sequence holds. 
Each homogeneous part of the sequence could be 
represented by matrices of $\QQ$-coefficients 
by any $\QQ$-basis, thus we have an algorithm to obtain $\eta$. 
\qed\end{pf}

\end{document}